\documentclass[preprint,12pt,number]{elsarticle}
\usepackage[margin=1in]{geometry}
\usepackage{amssymb}
\usepackage{amsmath}
\usepackage[unicode]{hyperref}
\hypersetup{colorlinks=true,
	citecolor=blue,
	linkcolor=blue}
\hypersetup{pdfauthor={Name}}
\usepackage{xcolor}
\usepackage{caption}
\usepackage{subcaption}
\usepackage{floatrow}
\newfloatcommand{capbtabbox}{table}[][\FBwidth]
\usepackage{blindtext}

\usepackage[skip=5pt plus1pt, indent=20pt]{parskip}

\usepackage[byauthor]{credits}
\credit{SS}{Methodology,Investigation,Software,Formal analysis,Validation,Writing -- original draft}
\credit{JM}{Supervision,Conceptualization,Methodology,Funding acquisition,Writing -- review \& editing}
\credit{KPKA}{Data curation,Investigation,Formal analysis,Writing -- review \& editing}
\credit{FL}{Data curation,Investigation,Formal analysis,Writing -- review \& editing}
\credit{SN}{Supervision,Conceptualization,Project administration,Funding acquisition,Writing -- review \& editing}
\credit{RD}{Project administration,Funding acquisition}

\usepackage{natbib}
\newcommand\blfootnote[1]{
    \begingroup
    \renewcommand\thefootnote{}\footnote{#1}
    \addtocounter{footnote}{-1}
    \endgroup
}
\newcommand{\RD}{\tilde{\rho}}
\newcommand{\fig}{\textcolor{blue}{Figure }}

\journal{Journal}

\begin{document}
\begin{frontmatter}


\title{A visco-plastic constitutive model for accurate densification and shape predictions in powder metallurgy hot isostatic pressing} 

\author[inst1]{Subrato Sarkar}
\author[inst1]{Jason R Mayeur}
\author[inst2]{KPK Ajjarapu}
\author[inst1]{Fred A List III}
\author[inst2]{Soumya Nag}
\author[inst1]{Ryan R Dehoff}
\affiliation[inst1]
    {organization={Manufacturing Science Division, 
        Oak Ridge National Laboratory},
            country={USA}}
\affiliation[inst2]{organization={Materials Science \& Technology Division, Oak Ridge National Laboratory},
            country={USA}}

\begin{abstract}
Powder metallurgy hot isostatic pressing (PM-HIP) is an advanced manufacturing process that produces near net shape parts with high material utilization and uniform microstructures. Despite being used frequently to produce small-scale components, the application of PM-HIP to large-scale components is limited due to inadequate understanding of its complex mechanisms that cause unpredictable post-HIP shape distortions. A computational model can provide necessary information about the intermediate and final stages of the HIP process that can help understand it better and make accurate predictions. Generally, two types of computational models are employed for PM-HIP simulations, namely, plastic and visco-plastic models. Between these, the plastic model is preferred due to its cheaper calibration approach requiring less experimental data. However, the plastic model sometimes produces incorrect predictions when slight variations of the HIP conditions are encountered in practical situations. Therefore, this work presents a visco-plastic model that addresses these limitations of the plastic model. A novel modified calibration approach is employed for the visco-plastic model that utilizes less experimental data than existing approaches. With the new approach, the data requirement is same for both plastic and visco-plastic models. This also enables a quantitative comparison of plastic and visco-plastic models, which have been only qualitatively compared in the past. When calibrated with the same experimental data, both the models are found to produce similar results. The calibrated visco-plastic model is applied to several complex geometries, and the predictions are found to be in good agreement with experimental observations.
\end{abstract}



\end{frontmatter}
\blfootnote{Notice: This manuscript has been authored by UT-Battelle, LLC, under contract DE-AC05-00OR22725 with the US Department of Energy (DOE). The US government retains and the publisher, by accepting the article for publication, acknowledges that the US government retains a nonexclusive, paid-up, irrevocable, worldwide license to publish or reproduce the published form of this manuscript, or allow others to do so, for US government purposes. DOE will provide public access to these results of federally sponsored research in accordance with the DOE Public Access Plan (https://www.energy.gov/doe-public-access-plan).}
\section{Introduction}
\label{intro}
Powder metallurgy hot isostatic pressing (PM-HIP) is an advanced manufacturing process that can produce parts with near net shape (NNS) and uniform microstructures while also reducing material wastage. A typical PM-HIP process involves: (1) Filling a can (canister/capsule) representing the desired component shape with metal powder, (2) Out gassing and sealing the can under vacuum, and (3) Subjecting the sealed can to a high pressure and temperature environment for several hours. An inert gas atmosphere is normally used to avoid oxidation at high temperatures of about 1100 $^{\circ}$C with pressures in the range of 100–120 MPa. PM-HIP is frequently used for producing smaller NNS parts with high precision and quality \cite{baccino2000high}, but its application to bigger parts is limited due to large, unpredictable post-HIP shrinkages and shape distortions. The shrinkages and distortions can be due to several factors such as, can stiffness, temperature gradients and powder characteristics \cite{Cassenti1980utcrept, atkinson2000fundamental, li1997constitutive}. Among these factors, the can stiffness and temperature gradients are directly affected by the can design \cite{svoboda1997simulation}. Therefore, it is crucial to understand the deformation mechanisms in the PM-HIP process for a suitable can design that can limit the shrinkages and shape distortions within acceptable ranges while also producing a fully densified part.

PM-HIP has been used in the past on various materials (such as metals, ceramics and composites) for a wide range of applications in aerospace, energy and medical fields \cite{lv2025review}. Despite extensive applications and industrial usage, the mechanistic understating of the PM-HIP process is still inadequate, leading to higher processing costs, longer developmental lead times and inability to effectively produce large sized or intricately shaped parts with good precision \cite{lv2025review}. The limited understanding of the PM-HIP process can be attributed to the presence of multiple deformation mechanisms and their complex coupled interactions. The deformation mechanisms include powder particle rearrangement, finite strain plasticity and creep under coupled thermo-mechanical loading conditions \cite{Cassenti1980utcrept}. Therefore, a mechanistic understanding of the PM-HIP process can help eliminate uncertainties associated with the prediction of post-HIP shrinkages and shape distortions. Using a mechanics-based computational model can significantly reduce developmental costs and lead times by reducing dependence on trial-and-error approaches \cite{lv2025review, atkinson2000fundamental}. 

The PM-HIP computational models in the past have been broadly classified into microscopic and macroscopic models \cite{li1997constitutive, atkinson2000fundamental}. The microscopic models study the powder densification response by analyzing a single powder particle and its surrounding, subject to pure hydrostatic stress under various densification mechanisms. In these models, powder densification equations are obtained by adding the effects of all physically occurring mechanisms at the microscale such as creep, grain boundary diffusion, volume diffusion and grain growth. As a result, these equations contain multiple parameters that are difficult to obtain through simple experiments such as uniaxial compression tests. Additionally, these densification equations are derived based on pure hydrostatic stresses and therefore are unable to predict powder response and shape distortion in complex geometries that experience appreciable deviatoric stresses during densification. Therefore, the application of microscopic models remains limited to isolated powder response and simple geometries with simple boundary conditions \cite{vanthe2016numerical}. 

In contrast, the macroscopic models consider metal powder as a continuum and study the overall densification behavior during the PM-HIP process \cite{Cassenti1980utcrept, abouaf1988finite}. In the macroscopic models, the classical theories of solid mechanics are modified to account for temperature-dependent effects in a porous continuum. The modifications mainly involve introducing additional parameters to characterize the temperature and density dependence of powder response during the PM-HIP process. Additionally, the macroscopic models are easier to calibrate using simple experiments and can be applied to a wide range of materials and geometries \cite{atkinson2000fundamental}. The macroscopic models have been shown to provide reasonably accurate predictions of powder densification and post-HIP shape distortions without needing ad hoc modifications like microscopic models \cite{atkinson2000fundamental}. Therefore, this work focusses on modeling the PM-HIP process using macroscopic models and applying them for predictions on complex geometries.

The majority of the macroscopic models used for PM-HIP simulations can be categorized into plastic, visco-plastic and combined (plastic + visco-plastic) models. The plastic models for PM-HIP are based on the classical \textit{J}$_2$ plasticity theory with modified yield function and hardening law to account for density and temperature dependence during the densification process \cite{shima1976plasticity, Cassenti1980utcrept}. Similarly, the visco-plastic models for PM-HIP use a modified form of Norton power law that accounts for the variations in the density and temperature during the densification process \cite{abouaf1988finite, svoboda1995determination, svoboda1996simulation}. A combined model has both plastic and visco-plastic models that are either active in parallel or serially \cite{van2017combined, wikman2000combined}. Among these, the plastic models are preferred due to their simpler calibration approach that requires less experimental data. However, the plastic models sometimes produce incorrect results when subject to slight variations in HIP conditions during densification, such as reduced applied pressure due to can stiffness or a variation in powder density at the beginning of densification. Therefore, a visco-plastic model is presented in this work that addresses these limitations of the plastic model. The main objectives of this work are to:

\begin{itemize}
    \item Present and validate the new modified calibration approach for the visco-plastic model that uses less experimental data than the existing approaches.
    \item Evaluate the performance of the calibrated visco-plastic model by comparing its densification behavior with the plastic model under slight variations in the HIP conditions.
    \item Apply the visco-plastic model on complex geometries and compare the predictions with the plastic model and experiments.
\end{itemize}

The remainder of the paper is organized into \hyperlink{Meth}{Methodology}, \hyperlink{RnD}{Results and Discussion} and \hyperlink{Con}{Conclusion} sections. The Methodology section discusses the constitutive models and the calibration process used for both plastic and visco-plastic models. The Results and Discussion section contrasts and compares the predictions of plastic and visco-plastic models for different geometries and conditions. Finally, the important observations and insights are summarized in the Conclusion section.

\section{Methodology}
\label{Meth}
In this section, firstly, the visco-plastic and plastic constitutive models used for PM-HIP modeling are described. Then, a detailed calibration approach for visco-plastic model is presented, which is a modified form of existing approaches. Finally, the validation of the calibrated visco-plastic model with experiments is carried out. Only a brief description of the plastic model is presented in this work because a detailed description of the formulation, computer implementation and calibration is presented in Mayeur \textit{et al.} \cite{mayeur2025MCCP}.
\subsection{Constitutive Models}
\label{contmodel}
This work uses a finite element (FE) based fully coupled thermo-mechanical model with finite deformations for simulating PM-HIP. The governing equations of the FE model based on momentum, energy and mass conservation respectively are:
\begin{equation}
    \nabla \cdot \boldsymbol{\sigma} = \mathbf{0}
\end{equation}
\begin{equation}
    \rho c_p \dot{T} = -\kappa \nabla^2 T
\end{equation}
\begin{equation}
    \dot{\rho} + \rho \, \nabla \cdot \mathbf{v} = 0
\end{equation}
where $ \boldsymbol{\sigma} $ is the Cauchy stress tensor, $ \rho $ is the mass density, $ c_p $ is the specific heat capacity at constant pressure, $ T $ is the temperature, $ \kappa $ is the thermal conductivity, and $ \mathbf{v} $ is the velocity vector. 

The constitutive models adopted in this work assume the metal powder compact (i.e. powder after filled into the can and compacted) as a porous continuum \cite{abouaf1988finite, svoboda1995determination, shima1976plasticity}. The powder porosity is represented by an internal variable called relative density ($\RD$). Relative density is the ratio of the porous material density ($ \rho $) to that of the fully dense material ($\rho_s$) expressed as $\RD=\frac{\rho}{\rho_s}$. The total strain rate is decomposed additively as,
\begin{equation}\label{eq:strainDecompP}
        \dot{\boldsymbol{\epsilon}} =\dot{\boldsymbol{\epsilon}}^e +\dot{\boldsymbol{\epsilon}}^\theta + 
        \begin{cases}
            \dot{\boldsymbol{\epsilon}}^{vp} & \text{Visco-plastic model} \\
            \dot{\boldsymbol{\epsilon}}^{p} & \text{Plastic model}
        \end{cases}
\end{equation}
where $\dot{\boldsymbol{\epsilon}}^e$ and $\dot{\boldsymbol{\epsilon}}^\theta$ stand for the elastic and thermal strains, while $\dot{\boldsymbol{\epsilon}}^{vp}$ and $\dot{\boldsymbol{\epsilon}}^{p}$ are the visco-plastic and plastic strains. Note that visco-plastic and plastic models are considered separate in this work, thus, the terms $  \dot{\boldsymbol{\epsilon}}^{vp} $ and $ \dot{\boldsymbol{\epsilon}}^{p} $ are active independently within the respective models. The model assumes isotropic and linear elastic response of the powder and the Cauchy stress tensor ($\boldsymbol{\sigma}$) is related to the elastic strain ($\boldsymbol{\epsilon^e}$) as,
\begin{equation}
    \boldsymbol{\sigma} = \lambda(\RD,T)\mathrm{tr}(\boldsymbol{\epsilon^e})\mathbf{I} + 2\mu(\RD,T)\boldsymbol{\epsilon^e}
\end{equation}
where, $ \mathbf{I} $ is the identity tensor and the Lame's parameters ($ \lambda(\RD,T) = E\nu/ \left[(1+\nu)(1-2\nu)\right] $ and $ \mu(\RD,T) = E/ \left[2(1+\nu)\right] $) depend on $\RD$ and $T$. The thermal strain in Eq. (\ref{eq:strainDecompP}) is defined as:
\begin{equation}\label{eq:strainThermal}
    \boldsymbol{\epsilon}^\theta = \alpha(\RD,T)\Delta T \, \mathbf{I}
\end{equation}
where $\alpha$ is the coefficient of thermal expansion (CTE) dependent on $\RD$ and $T$. 

%
\subsubsection{Visco-plastic model}
This work uses a modified form of visco-plastic model proposed by Abouaf \textit{et al.} \cite{abouaf1988finite} that accounts for the dependence of the model on relative density ($ \RD $). The visco-plastic strain rate ($\dot{\boldsymbol{\epsilon}}^{vp}$) is defined as,
\begin{equation}\label{eq:strainVP}
    \dot{\boldsymbol{\epsilon}}^{vp} = 
        A(T)\sigma_{eqv}^{{N(T)-1}} \left(\frac{3}{2} c(\RD)\mathbf{s} + f(\RD)I_1\mathbf{I}\right)
\end{equation}
where, $A(T)$ and $N(T)$ are the temperature dependent material constants of power-law type visco-plastic behavior. $c(\RD)$ and $f(\RD)$ are the relative density dependent weighing functions that determine the contributions of the deviatoric ($ \mathbf{s} $) and hydrostatic ($I_1=\mathrm{tr}(\boldsymbol{\sigma})$) stress on the strain rate during the densification process. The equivalent stress in Eq. (\ref{eq:strainVP}) is defined as a function of $\RD$,
\begin{equation}\label{eq:eqStressCr}
    \sigma_{eqv}(\RD) = \left(3c(\RD)J_2 + f(\RD)I_1^2\right)^{\frac{1}{2}}
\end{equation}
where, $J_2 = \frac{1}{2}\mathbf{s}:\mathbf{s}$, is the second invariant of the deviatoric stress tensor.

The visco-plastic model is implemented in Abaqus using the user-defined CREEP subroutine, which is suitable for models with strain rate dependence on both deviatoric and hydrostatic stress \cite{dunne2005introduction, vanthe2016numerical}. The CREEP subroutine requires the definition of equivalent swelling ($\Bar{\epsilon}^{sw}$) and creep ($\Bar{\epsilon}^{cr}$) strain increments and their derivatives to effectively use the implicit time integration scheme \cite{abaqus2024}. The incremental total visco-plastic strain in the subroutine is defined as,
\begin{equation}\label{eq:AbqDeltaStain}
    \Delta\boldsymbol{\epsilon}^{vp} = 
        \frac{1}{3}\Delta\Bar{\epsilon}^{sw}\mathbf{I} +
        \Delta\Bar{\epsilon}^{cr}\mathbf{n}
\end{equation}
where, $\mathbf{n}=\frac{\partial{\sigma_{eq}}}{\partial{\boldsymbol{\sigma}}}$, is called the deviatoric stress potential. Note that $\sigma_{eq}$ in deviatoric stress potential is the von-Mises equivalent stress which is different from  $\sigma_{eqv}$ defined in Eq. (\ref{eq:eqStressCr}). The incremental swelling and creep strains in Eq. (\ref{eq:AbqDeltaStain}) are defined using Eq. (\ref{eq:strainVP}) as,
\begin{equation}\label{eq:AbqDeltaStainSw}
    \Delta\Bar{\epsilon}^{sw} = \left( 3f(\RD)I_1 A(T) {\sigma}_{eqv}^{N(T)-1} \right) \Delta{t}
\end{equation}
\begin{equation}\label{eq:AbqDeltaStainCr}
    \Delta\Bar{\epsilon}^{cr} = \left( c(\RD)\sigma_{eq} A(T) {\sigma}_{eqv}^{N(T)-1} \right) \Delta{t}
\end{equation}
where, $ \Delta{t} $ is time increment. The derivatives of $\Delta\Bar{\epsilon}^{sw}$ and $\Delta\Bar{\epsilon}^{cr}$ that are needed for the implicit time integration are defined in \ref{app_deriv}. The relative density ($ \RD $) evolution is governed by visco-plastic volumetric strain expressed as,
\begin{equation}\label{eq:RD_Evol}
    \RD = \RD_o e^{\boldsymbol{-\epsilon}_{vol}^{vp}}
\end{equation}
where, $ \RD_o $ is the initial relative density of the powder compact before HIPing and the visco-plastic volumetric strain ($\boldsymbol{\epsilon}_{vol}^{vp}=\mathrm{tr}(\boldsymbol{\epsilon}^{vp})$). The variable $ \RD $ is not directly available in the CREEP subroutine. Hence, it is defined and implemented using field variables through the user-defined field (USDFLD) subroutine available in Abaqus. The Eqs. (\ref{eq:AbqDeltaStain})-(\ref{eq:RD_Evol}) clearly indicate that determination of parameters $A(T)$, $N(T)$, $c(\RD)$ and $f(\RD)$ is essential for a fully defined visco-plastic PM-HIP model in Abaqus. These parameters are determined using experimental data through a calibration process described later in this work.
%
\subsubsection{Plastic model}
The plastic model used in this work uses a modified form of  $J_2$ plasticity theory that accounts for the dependence on hydrostatic stress component ($I_1$) and relative density ($\RD$) during the PM-HIP process \cite{shima1976plasticity}. An complete description of the formulation and implementation of the PM-HIP plastic model using Modified Cam-Clay Plasticity (MCCP) is provided in Mayeur \textit{et al.} \cite{mayeur2025MCCP}. In the MCCP model, the evolution of an elliptical yield surface during the PM-HIP is defined in $p-q$ (pressure-deviatoric) stress space as, 
\begin{equation}\label{eq:camClayYield}
    \frac{p^2}{a^2} + \frac{q^2}{(Ma)^2} = 1
\end{equation}
where, $p = -I_1/3$ is hydrostatic pressure and $q=\sqrt{3 J_2}$ is the deviatoric stress. The parameters $a$ and $Ma$ in Eq. (\ref{eq:camClayYield}) are the semi-axis lengths along $p$ (pressure) and $q$ (deviatoric) axes of the elliplic yield surface, respectively. The temperature dependence of the yield surface evolution is defined through a hardening response in terms of pressure ($p$) and plastic volumetric strain ($\epsilon^p_{vol}$). For implementation in Abaqus, the parameter $M(\RD,T)$ is defined using *CLAY PLASTICITY with an associative hardening flow rule (in terms of $p-\epsilon^p_{vol}$ response) defined using *CLAY HARDENING material definitions. Consequently, the calibration of plastic model involves determination of $M(\RD,T)$ and $p-\epsilon^p_{vol}$ response using experimental observations \cite{mayeur2025MCCP}. Similar to the visco-plastic model, the dependence of different parameters on $\RD$ is defined through the user-defined field (USDFLD) subroutine in Abaqus.

\subsection{Visco-plastic Model Calibration}
%
\subsubsection{Existing and current calibration approaches}
\begin{figure}[H]
    \centering
    \includegraphics[width=0.75\textwidth]{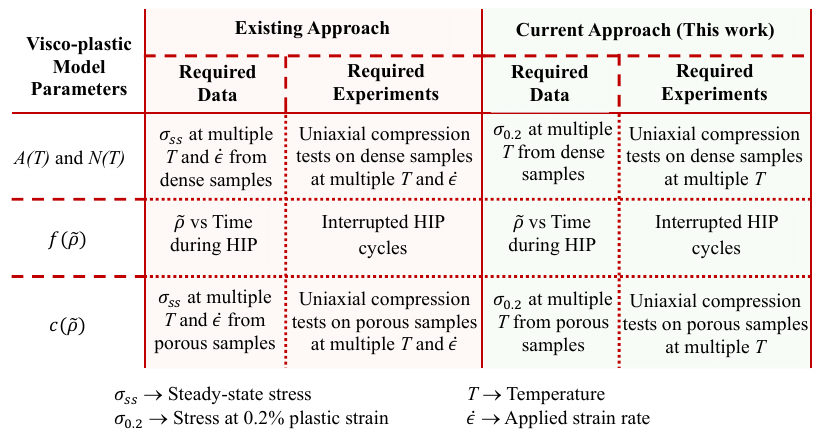}
    \caption{A comparison of the current calibration approach adopted in this work with existing approach for visco-plastic PM-HIP model parameters. The current approach uses $\sigma_{0.2}$ instead of $\sigma_{ss}$ and does not need additional data at different $\dot{\epsilon}$ for $A(T)$, $N(T)$ and $c(\RD)$.}
    \label{fig:calib_comp}
\end{figure}

As outlined in the previous section, the visco-plastic PM-HIP model implemented in Abaqus using the CREEP subroutine requires determination of four parameters, i.e. $A(T)$, $N(T)$, $c(\RD)$ and $f(\RD)$, that are either dependent on $T$ or $\RD$. The existing calibration approach used for visco-plastic PM-HIP models has been discussed in detail by Svoboda \textit{et al.} \cite{svoboda1995determination, svoboda1996simulation} and applied by various others successfully \cite{van2017combined, abouaf1988finite, cho2001densification}. The experimental data and associated experiments required in the existing calibration approach are outlined in \fig \ref{fig:calib_comp}. For calibrating the transient densification behavior ($c(\RD)$) during HIP, density measurements of partially consolidated samples (interrupted HIP cycles) are needed. While the uniaxial isothermal compression tests at different applied temperatures and strain rates are needed to establish the flow behavior ($A(T)$, $N(T)$ and $f(\RD)$) of metal powder during HIP. It is pointed out that determining flow behavior is expensive as it requires testing at different strain rates. The current calibration approach adopted in this work is slightly different and the major differences are:

\begin{itemize}
    \item Using uniaxial compression test data at a single strain rate instead of multiple strain rates for determining $A(T)$, $N(T)$ and $c(\RD)$, thus, reducing experimentation costs.
    \item Using 0.2\% stress ($\sigma_{0.2}$) instead of steady-state stress ($\sigma_{ss}$) from the uniaxial compression tests for determining $A(T)$, $N(T)$ and $c(\RD)$, thus, establishing equivalence with the plastic model.
\end{itemize}

\begin{figure}[h]
    \centering
    \includegraphics[width=0.75\textwidth]{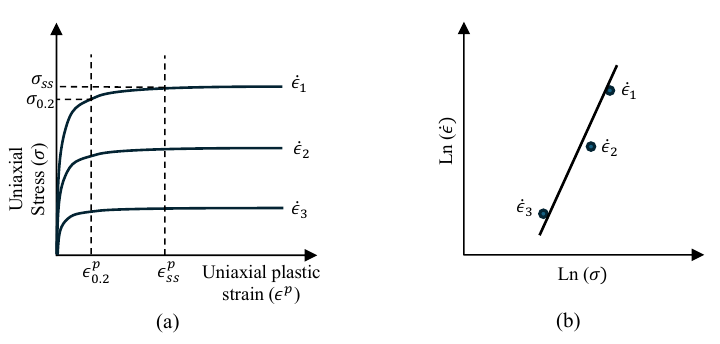}
    \caption{A schematic showing the typical visco-plastic stress-strain response 
    in a uniaxial compression test through \textbf{(a)} stress-strain plot at different strain rates and \textbf{(b)} log-log plot of strain rate vs stress.}
    \label{fig:sche1}
\end{figure}

The differences can be understood better from the schematic plots shown in \fig \ref{fig:sche1}. \fig \ref{fig:sche1}a shows the stress-strain response of a typical visco-plastic behavior observed in uniaxial compression tests. Note that 0.2\% ($\sigma_{0.2}$) and steady-state ($\sigma_{ss}$) stresses are obtained at 0.2\% plastic strain ($\epsilon_{0.2}^p$) and steady-state plastic strain ($\epsilon_{ss}^p$). The steady-state behavior is observed when changes in the plastic strain do not cause significant change in the stress. In a visco-plastic material, different applied strain rates ($\dot{\epsilon}$) lead to different stress-strain responses, as shown on a log-log plot in \fig \ref{fig:sche1}b. The points on the log-log plot are used to fit a straight line for determining the visco-plastic model parameters. Normally, $\sigma_{ss}$ is used in the log-log plots, however, $\sigma_{0.2}$ is used in this work instead. The $\sigma_{0.2}$ is found suitable because, at higher temperatures ($\geq800^o$C), $\sigma_{0.2}$ and $\sigma_{ss}$ are equivalent in compression tests \cite{vanthe2016numerical}. Additionally, $\sigma_{0.2}$ is considered as the yield stress in compression and used to calibrate the yielding behavior in visco-plastic model. The reasoning and conditions for using a single strain rate compression test data for calibration are described later in this section. SS 316L is used as the model material for calibration of the visco-plastic PM-HIP model mainly due to the availability of all required experimental data \cite{vanthe2016numerical, kohar2019new}.

%
\subsubsection{Calibration steps and validation}

\begin{figure}[h]
    \centering
    \begin{subfigure}[b]{0.375\textwidth}
        \centering
        \includegraphics[width=\textwidth]{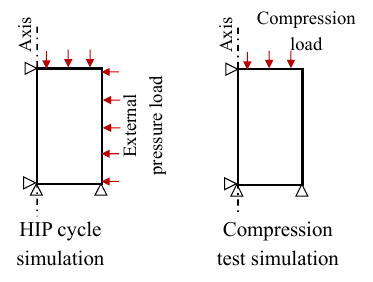}
        \caption{Axisymmetric single element models}
        \label{}
    \end{subfigure}
    \hspace{1cm}
    \begin{subfigure}[b]{0.375\textwidth}
        \centering
        \includegraphics[width=\textwidth]{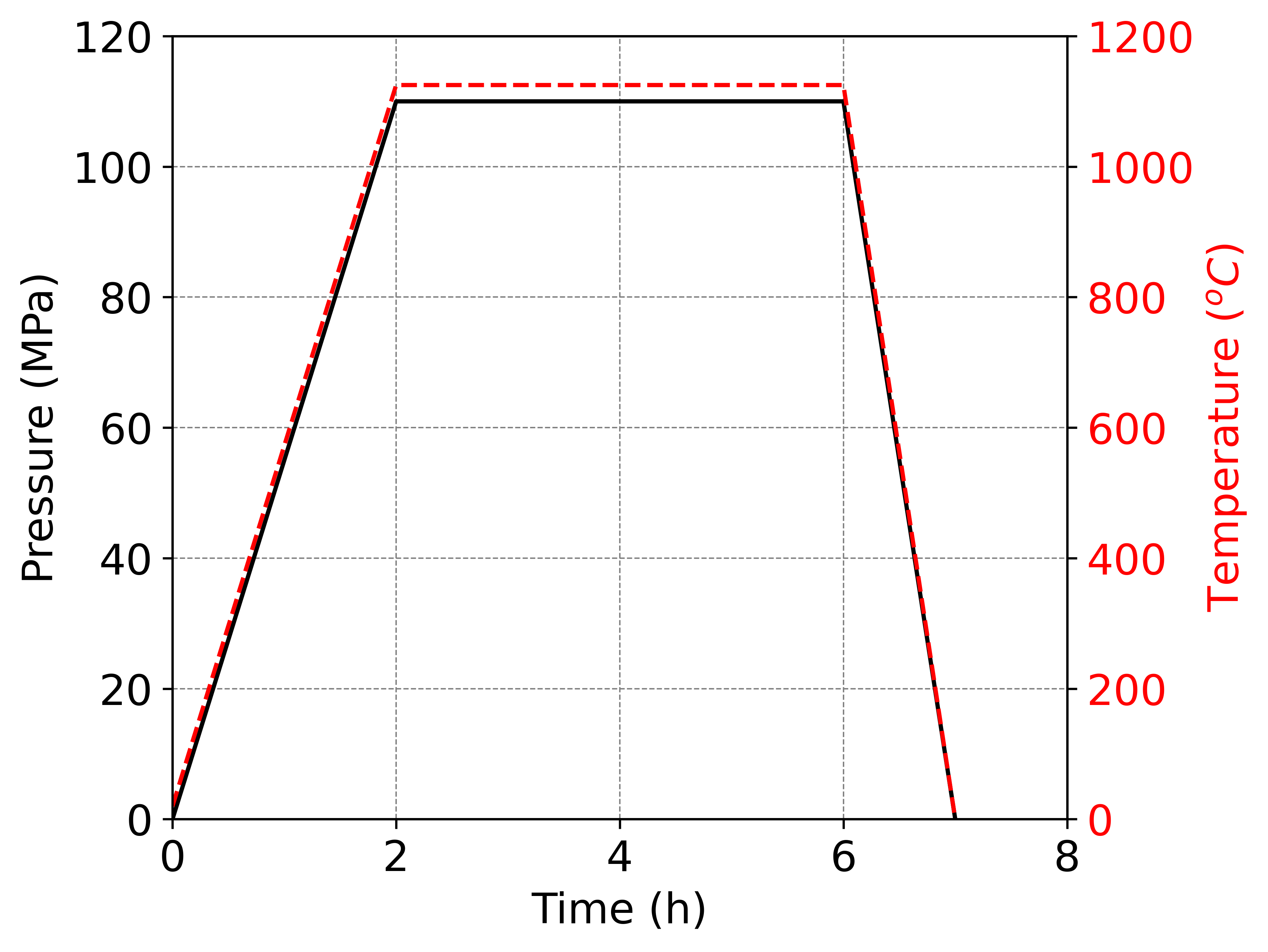}
        \caption{HIP cycle pressure and temperature}
        \label{}
    \end{subfigure}
    \caption{ \textbf{(a)} The single element models used for calibration and \textbf{(b)} the HIP cycle adopted from Nguyen \textit{et al.} \cite{vanthe2016numerical}. }
    \label{fig:singElem_HIPCyc}
\end{figure}

The steps involved in the current calibration approach are similar to the existing approach, except for the differences mentioned in the previous subsection. The main steps in the current approach are outlined below:

\begin{enumerate}
    \item Calibrate $A(T)$ and $N(T)$ using yield (0.2\%) stress obtained from compression tests of the fully dense samples at different temperatures
    \item Calibrate $f(\RD)$ using the densification rate data obtained from density measurements of samples from interrupted HIP cycles
    \item Calibrate $c(\RD)$ using yield (0.2\%) stress obtained from compression tests of the partially consolidated (porous) samples at different temperatures
\end{enumerate}
The simulations for the above-mentioned calibration steps are carried out on a single element model (\fig \ref{fig:singElem_HIPCyc}a) in Abaqus having dimensions similar to the compression test specimen reported by Nguyen \textit{et al.} \cite{van2017combined}. The HIP cycle temperature and pressure used for the calibration are adopted from Nguyen \textit{et al.} \cite{van2017combined} and shown in \fig \ref{fig:singElem_HIPCyc}b. The thermoelastic properties of the SS316 powder dependent on $\RD$ and $T$ are adopted from literature \cite{vanthe2016numerical,kohar2019new}. It is pointed out that the calibration is performed only for the powder response, while, the can material is assumed to have rate-independent plastic response, for which material properties are adopted from literature \cite{vanthe2016numerical,kohar2019new,mayeur2025MCCP}. Each step in the calibration process is described in the following.

\noindent\textbf{Step 1: $A(T)$ and $N(T)$ calibration}

\begin{figure}[h] 
    \centering
    \begin{subfigure}[b]{0.36\textwidth}
        \centering
        \includegraphics[width=\textwidth]{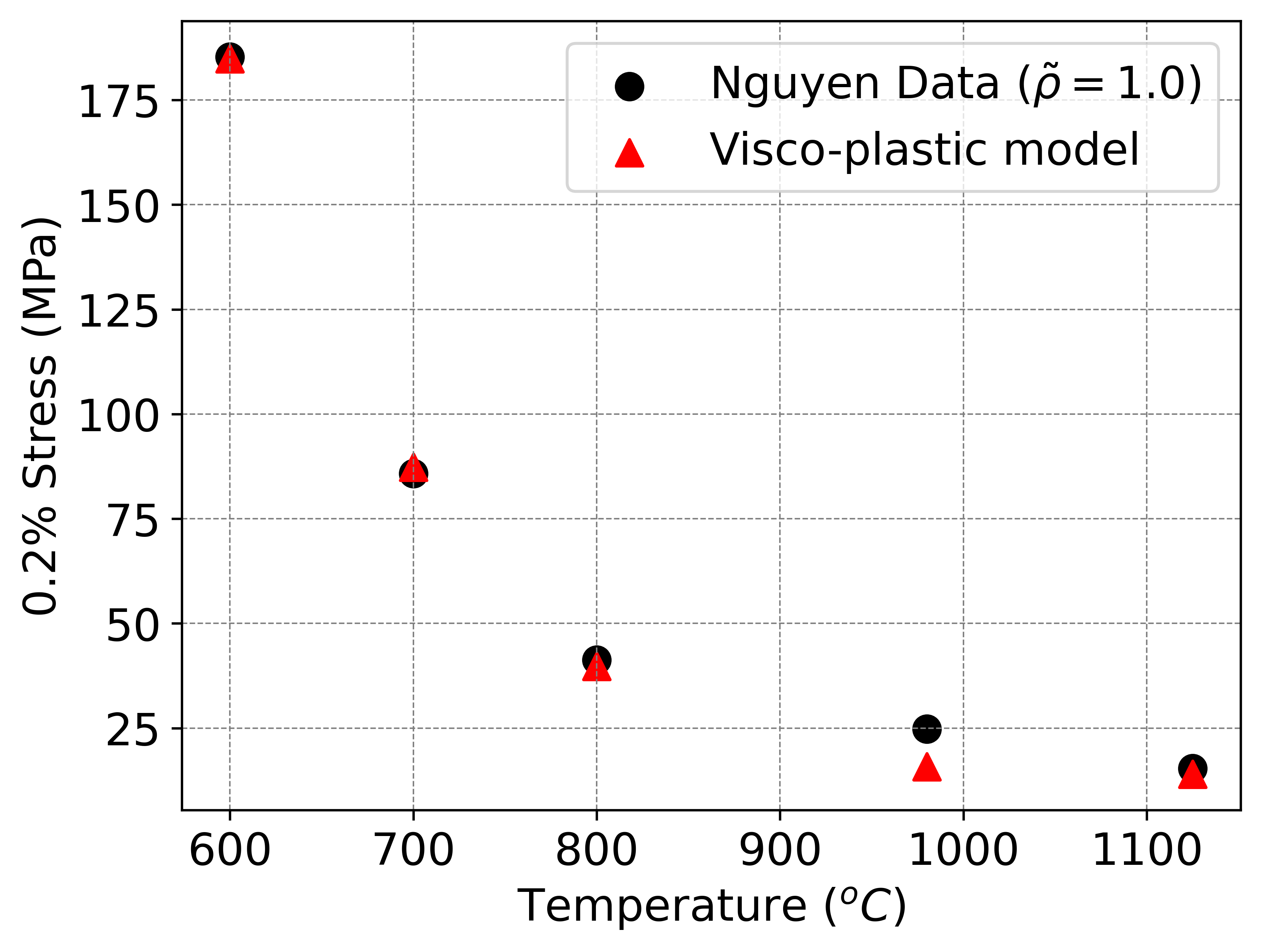}
        \caption{Yield stress validation (Dense samples)}
        \label{}
    \end{subfigure}
    \hspace{1cm}
    \begin{subfigure}[b]{0.375\textwidth}
        \centering
        \includegraphics[width=\textwidth]{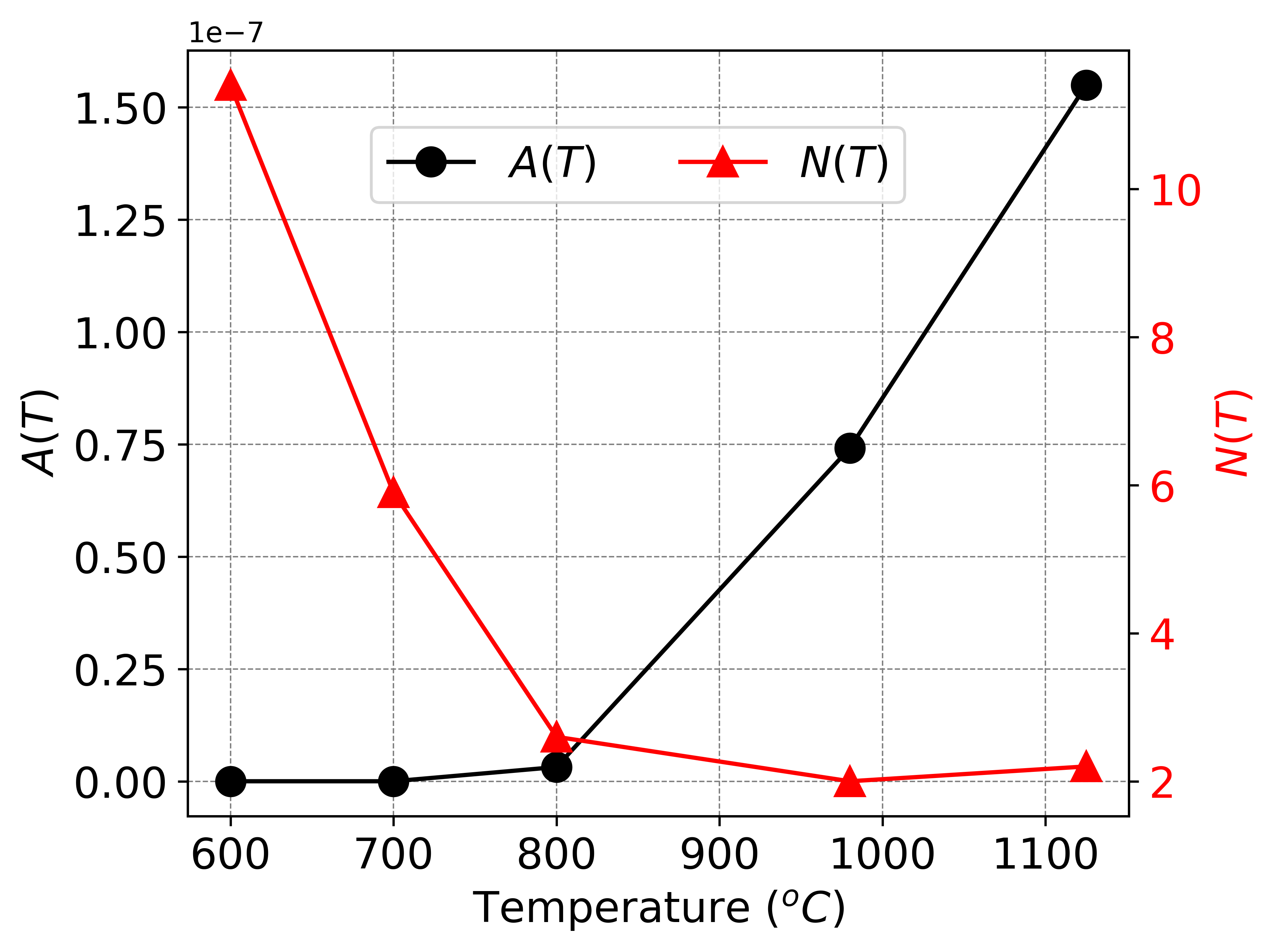}
        \caption{Calibrated $A(T)$ and $N(T)$}
        \label{}
    \end{subfigure}
    \caption{The results of $A(T)$ and $N(T)$ calibration showing \textbf{(a)} validation of yield (0.2\%) stress obtained from the visco-plastic model with the reference data (Nguyen \textit{et al.} \cite{vanthe2016numerical}) at different temperatures for fully dense samples and \textbf{(b)} final calibrated values of $A(T)$ and $N(T)$ across different temperatures.}
    \label{fig:ANcalibPlots}
\end{figure}

For calibration of $A(T)$ and $N(T)$, the single element model is considered fully dense ($\RD=1.0$) and is subject to uniaxial compression. In a fully dense material, the parameters $c(\RD)=1.0$ and $f(\RD)=0.0$ transform Eqs. (\ref{eq:strainVP}) and (\ref{eq:eqStressCr}) to,

\begin{equation}\label{eq:strainVP_ANcal}
    \dot{\boldsymbol{\epsilon}}^{vp} = 
        \frac{3}{2} A(T)\sigma_{eqv}^{{N(T)-1}} \mathbf{s}
\end{equation}
\begin{equation}\label{eq:eqStressCr_ANcal}
    \sigma_{eqv}(\RD) = (3J_2)^{\frac{1}{2}} = |\sigma_1|
\end{equation}

\noindent where $\sigma_1$ is the uniaxial compression stress. Using Eqs. (\ref{eq:strainVP_ANcal}) and (\ref{eq:eqStressCr_ANcal}) in the single element simulations, the unknowns $A(T)$ and $N(T)$ are determined by reducing the error between yield (0.2\%) stress obtained from simulations and the experimental compression yield stress reported by Nguyen \textit{et al.} \cite{vanthe2016numerical}. The error reduction is achieved by optimization of $A(T)$ and $N(T)$ resulting in close agreement of 0.2\% stress from simulations with the experiments as shown in \fig \ref{fig:ANcalibPlots}a. The 0.2$\%$ stress values from simulations (visco-plastic model) shown in \fig \ref{fig:ANcalibPlots}a are obtained by using calibrated $A(T)$ and $N(T)$ combinations shown in \fig \ref{fig:ANcalibPlots}b at different temperatures. It must be noted that compression test data at a single strain rate and already known general creep behavior for SS 316L are used to arrive at the initial guess of $A(T)$ and $N(T)$ for optimization. In the absence of known behavior, compression tests at atleast two strain rates might become necessary. In this work, the values of $A(T)$ and $N(T)$ are interpolated linearly between the calibrated temperature points shown in \fig \ref{fig:ANcalibPlots}b without fitting any smooth curves. For $T<600^o$C, $A(T)$ and $N(T)$ are assumed to be constant and equal to the values at $T=600^o$C. 

\noindent\textbf{Step 2: $f(\RD)$ calibration}

\begin{figure}[H] 
    \centering
    \begin{subfigure}[b]{0.36\textwidth}
        \centering
        \includegraphics[width=\textwidth]{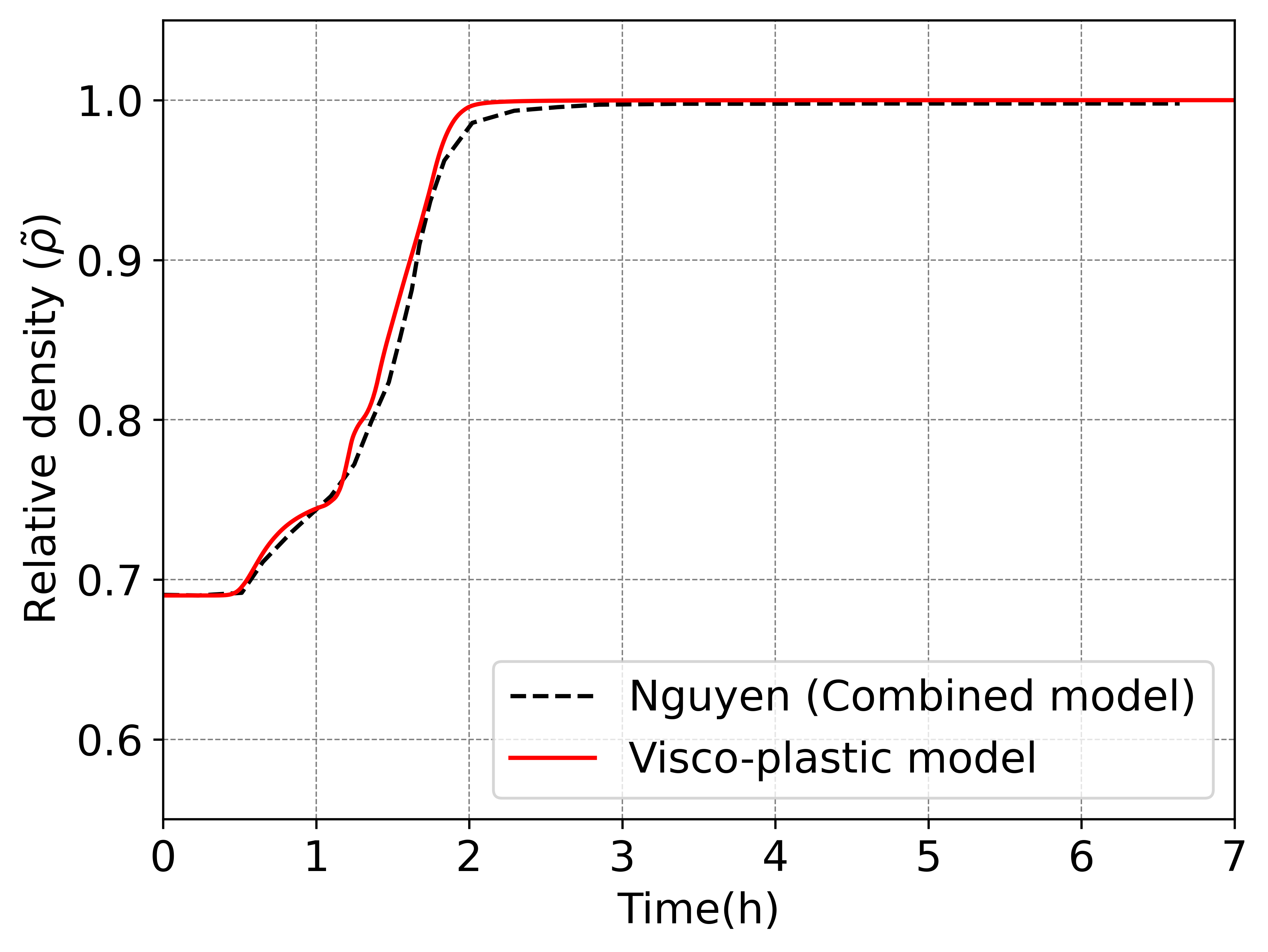}
        \caption{Densification response validation}
        \label{}
    \end{subfigure}
    \hspace{1cm}
    \begin{subfigure}[b]{0.36\textwidth}
        \centering
        \includegraphics[width=\textwidth]{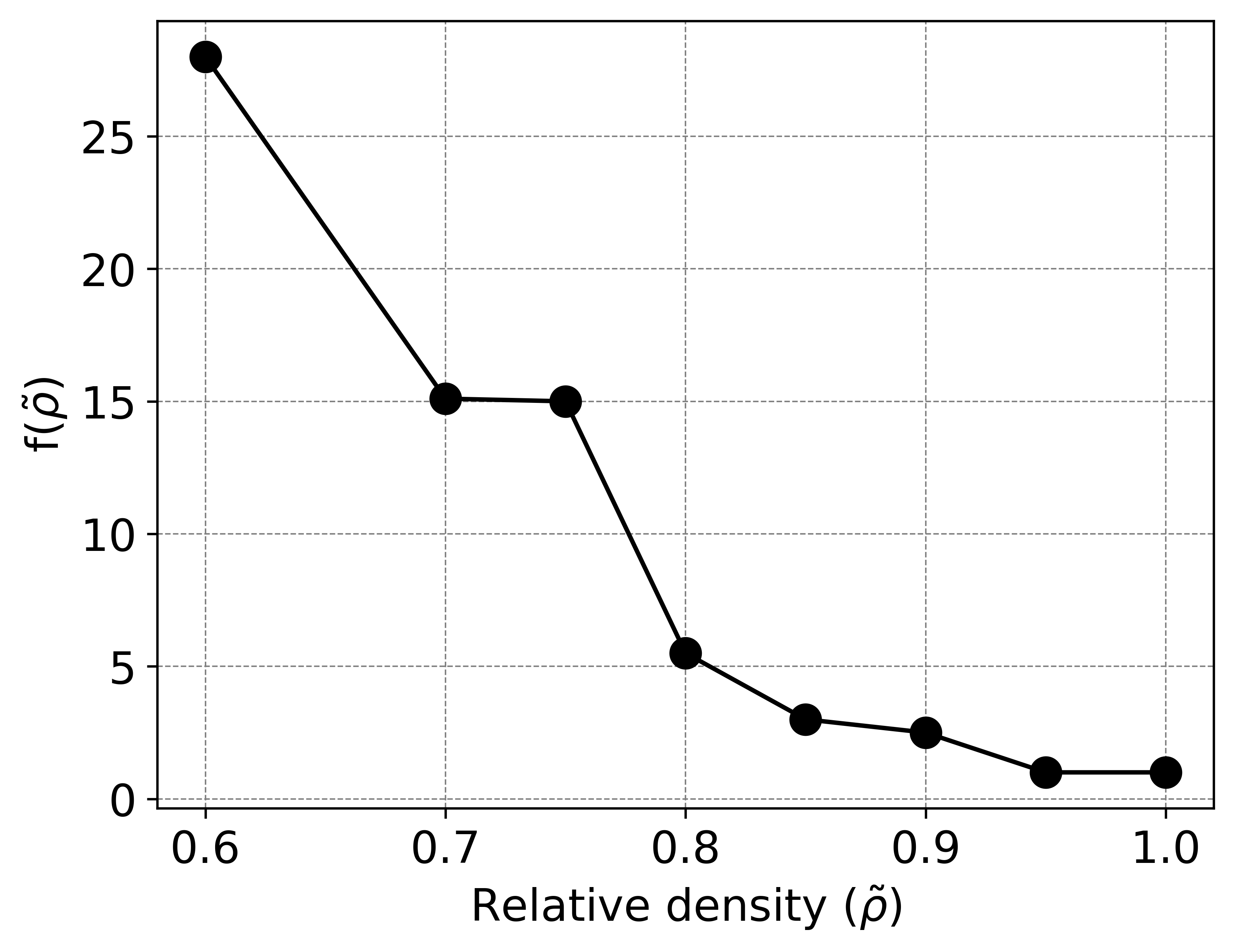}
        \caption{Calibrated $f(\RD)$}
        \label{}
    \end{subfigure}
    \caption{The results of $f(\RD)$ calibration showing \textbf{(a)} validation of densification response obtained from the visco-plastic model with the reference data (Nguyen \textit{et al.} \cite{vanthe2016numerical}) and \textbf{(b)} final calibrated values of $c(\RD)$ across the entire range of observed relative density.}
    \label{fig:fcalibPlots}
\end{figure}

With $A(T)$ and $N(T)$ fixed from the previous step, single element densification simulations under HIP cycle pressure and temperature (See \fig \ref{fig:singElem_HIPCyc}b) are carried out to determine $f(\RD)$. Due to the simple geometry, the single element model is under pure hydrostatic stress when subjected to HIP cycle pressures. Under pure hydrostatic stress, the deviatoric components in Eqs. (\ref{eq:strainVP}) and (\ref{eq:eqStressCr}) vanish to yield the following expressions,
\begin{equation}\label{eq:strainVP_fcal}
    \dot{\boldsymbol{\epsilon}}^{vp} = 
        A(T)\sigma_{eqv}^{{N(T)-1}} f(\RD)I_1\mathbf{I}
\end{equation}
\begin{equation}\label{eq:eqStressCr_fcal}
    \sigma_{eqv}(\RD) = \left(f(\RD)I_1^2\right)^{\frac{1}{2}}
\end{equation}
Using Eqs. (\ref{eq:RD_Evol}), (\ref{eq:strainVP_fcal}) and (\ref{eq:eqStressCr_fcal}), an expression of $f(\RD)$ is obtained as,
\begin{equation}\label{eq:fcal_eq}
    f(\RD) = \left(\frac{\dot{\RD}}{\RD A(T)(I_1)^{N(T)}}\right)^{\frac{2}{N(T)+1}}
\end{equation}
In Eq. (\ref{eq:fcal_eq}), $A(T)$ and $N(T)$ are known quantities and the unknown $f(\RD)$ is determined using the densification response obtained experimentally from density measurements of partially consolidated HIP samples \cite{svoboda1995determination}. However, in the absence of in-house experimental data, the densification curve of the combined model in Nguyen \textit{et al.} \cite{van2017combined} is used as a reference. \fig \ref{fig:fcalibPlots}a shows the densification response predicted by the single element model, and it agrees well with the combined model of Nguyen \textit{et al.} \cite{van2017combined}. The values of  $f(\RD)$ used to get the matching densification response are shown in \fig \ref{fig:fcalibPlots}b. Note that similar to $A(T)$ and $N(T)$, a linear interpolation is used to compute values of $f(\RD)$ at intermediate temperatures.

\noindent\textbf{Step 3: $c(\RD)$ calibration}

Parameter $c(\RD)$ is determined using uniaxial compression simulations on the single element model representing partially consolidated (porous) samples. The known values of $A(T)$, $N(T)$ and $f(\RD)$ from previous steps are used, and the yield (0.2\%) stresses from simulations are compared with experimental compression yield stress reported by Nguyen \textit{et al.} \cite{vanthe2016numerical}. Under uniaxial compression stress ($\sigma_1$), the stress invariants can be expressed as $I_1=\sigma_1$ and $\sqrt{3J_2}=|\sigma_1|$ transforming the Eqs. (\ref{eq:strainVP_fcal}) and (\ref{eq:eqStressCr_fcal}) into,

\begin{equation}\label{eq:strainVP_ccal}
    {\dot{\epsilon}_1}^{vp} = 
        A(T)\sigma_1^{{N(T)}} \left[c(\RD) + f(\RD)\right]^\frac{N(T-1)}{2}
\end{equation}
\begin{equation}\label{eq:eqStressCr_ccal}
    \sigma_{eqv} = |\sigma_1| \sqrt{c(\RD) + f(\RD)}
\end{equation}

\begin{figure}[h] 
    \centering
    \begin{subfigure}[b]{0.36\textwidth}
        \centering
        \includegraphics[width=\textwidth]{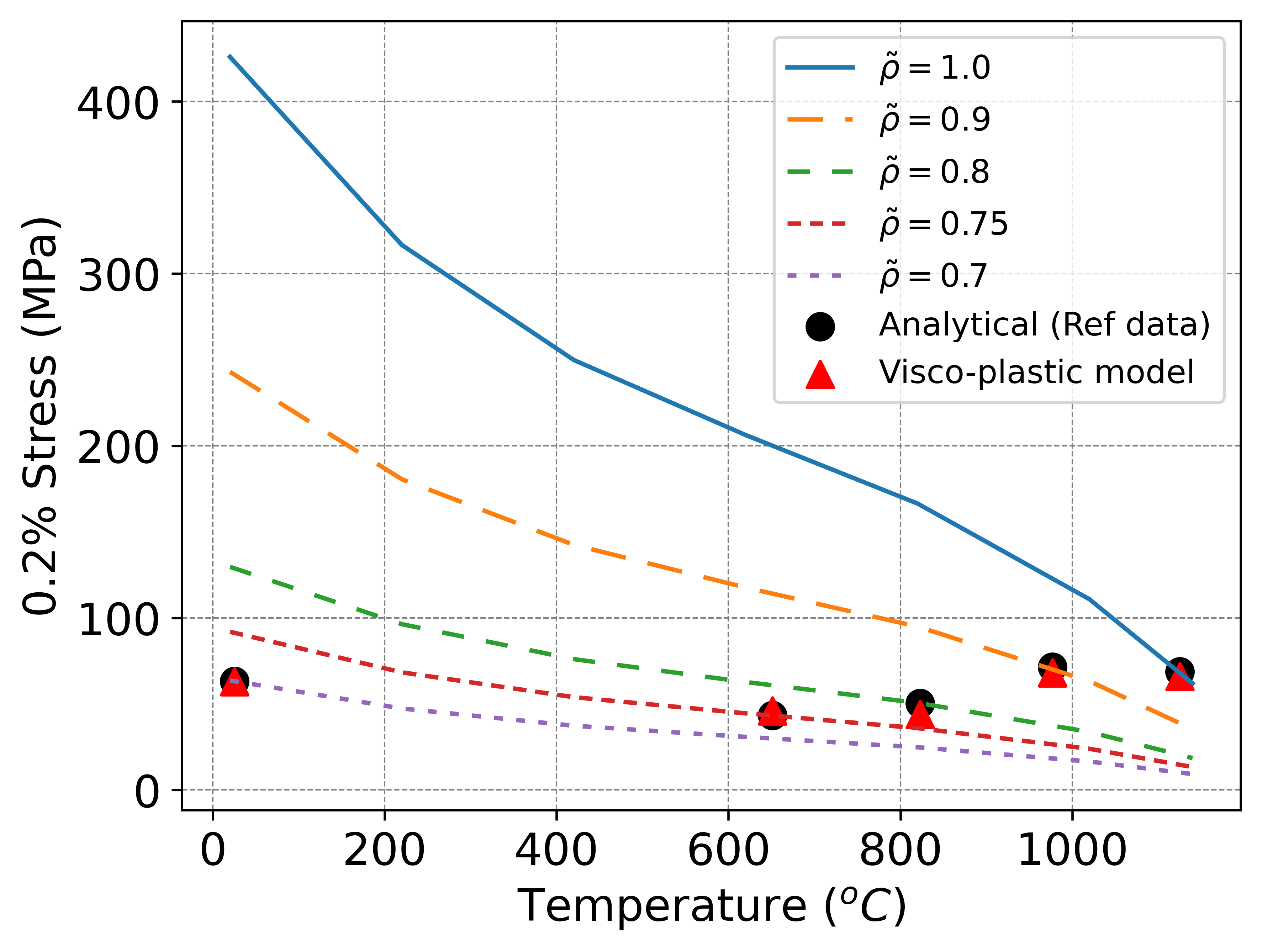}
        \caption{Yield stress validation (Porous samples)}
        \label{}
    \end{subfigure}
    \hspace{1cm}
    \begin{subfigure}[b]{0.36\textwidth}
        \centering
        \includegraphics[width=\textwidth]{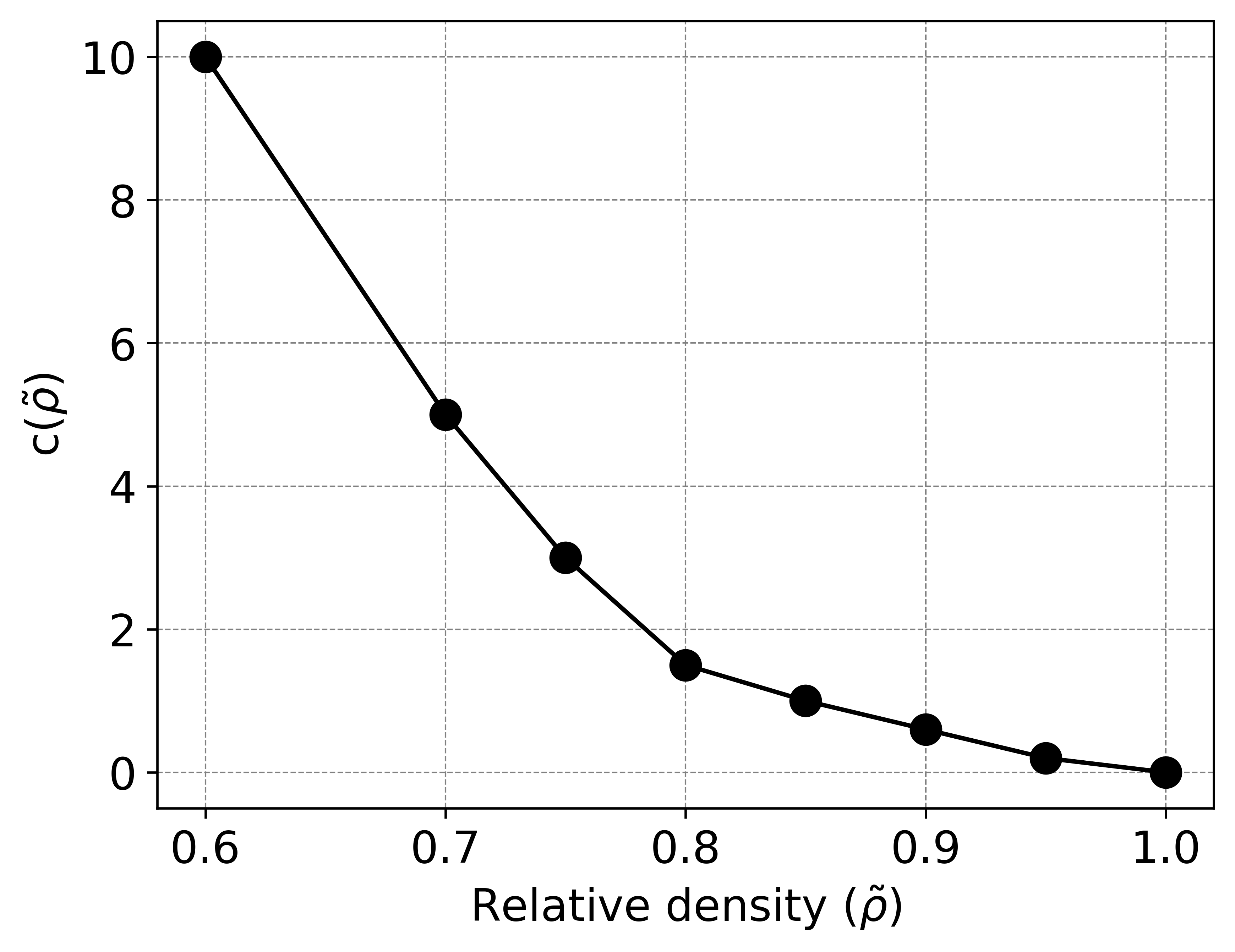}
        \caption{Calibrated $c(\RD)$}
        \label{}
    \end{subfigure}
    \caption{The results of $c(\RD)$ calibration showing \textbf{(a)} validation of yield (0.2\%) stress response obtained from the visco-plastic model for partially consolidated (porous) samples with the reference data (Nguyen \textit{et al.} \cite{vanthe2016numerical}) and \textbf{(b)} final calibrated values of $f(\RD)$ across the entire range of observed relative density.}
    \label{fig:ccalibPlots}
\end{figure}

\noindent Using Eq. (\ref{eq:strainVP_ccal}), an expression of $c(\RD)$ is obtained as,

\begin{equation}\label{eq:ccal_eq}
    c(\RD) = \left(\frac{\dot{\epsilon}_1}{A(T) {\sigma_1}^{N(T)}}\right)^{\frac{2}{N(T)-1}}-f(\RD)
\end{equation}

\noindent where the subscript 1 represents uniaxial values. Similar to step 1, $c(\RD)$ is determined through the minimization of 0.2\% stress at different temperature and relative density combinations encountered during the HIP cycle (\fig \ref{fig:singElem_HIPCyc}b). \fig \ref{fig:ccalibPlots}a shows analytical curves fitted to the experimental data at different $\RD$ in Nguyen \textit{et al.} \cite{vanthe2016numerical} and the black dots are the sample points corresponding to temperature ($T$) and relative density ($\RD$) combinations encountered during the HIP cycle. An iterative process is adopted to determine the $c(\RD)$ at various $T$ and $\RD$ combinations by minimizing the error between simulated and experimental yield stress. The initial values of $c(\RD)$ at the start of the minimization process are determined by substituting 0.2\% stress and the strain rate in Eq. (\ref{eq:ccal_eq}) from the uniaxial compression tests on porous samples at different temperatures. The simulations are repeated until the error between 0.2\% stress is minimized, as shown with red triangles in \fig \ref{fig:ccalibPlots}a. The final $c(\RD)$ values obtained after the minimization are shown in \fig \ref{fig:ccalibPlots}b.

The determination of $c(\RD)$ using Eq. (\ref{eq:ccal_eq}) typically requires uniaxial compression tests at multiple strain rates, as demonstrated in previous works \cite{van2017combined, svoboda1995determination}. However, it is found that the multiple strain rates data is only useful for computing the initial guess of parameter $c(\RD)$ and the final calibrated values are significantly different from the initial guess. Therefore, compression tests at multiple strain rates may not be required, if a reasonably good guess can be made using only the single strain rate data and by studying general creep behavior of the material of interest from other sources. However, multiple (or at least two) strain rate tests might become necessary when sufficient data us unavailable.

\section{Results and Discussion}
\label{RnD}
A comparison of the plastic and the visco-plastic model is presented in this section through simulating HIP of several different geometries. The results are compared with either experimental observations or reference data available in the literature. The final calibrated parameters for both the models are used for simulations in this section. The calibrated parameters for the plastic (modified cam-clay) model are listed in Mayeur \textit{et al.} \cite{mayeur2025MCCP}.
\subsection{Model generalization under different initial relative density }
\begin{figure}[H] 
    \centering
    \begin{subfigure}[b]{0.36\textwidth}
        \centering
        \includegraphics[width=\textwidth]{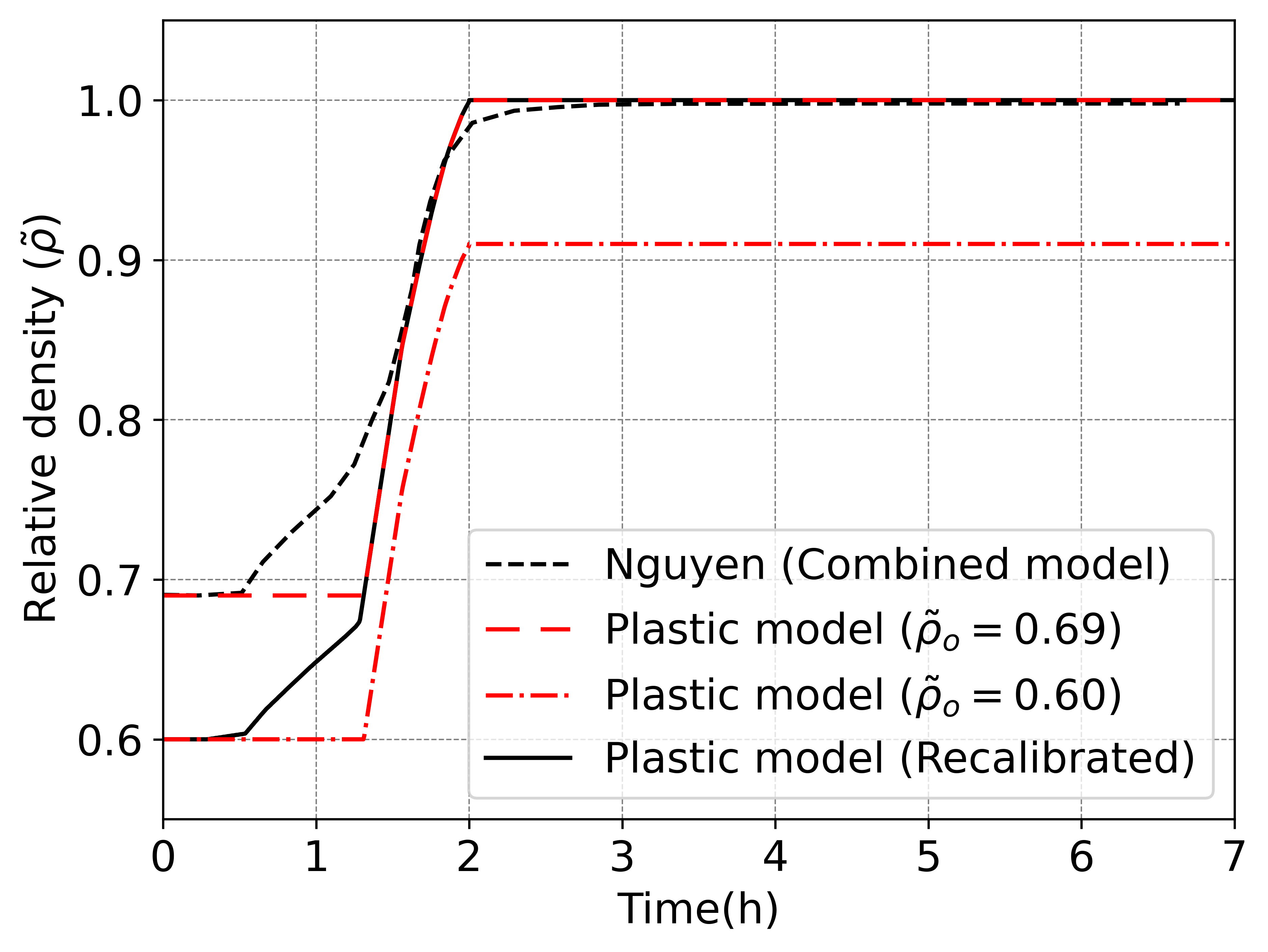}
        \caption{Plastic model}
        \label{}
    \end{subfigure}
    \hspace{1cm}
    \begin{subfigure}[b]{0.36\textwidth}
        \centering
        \includegraphics[width=\textwidth]{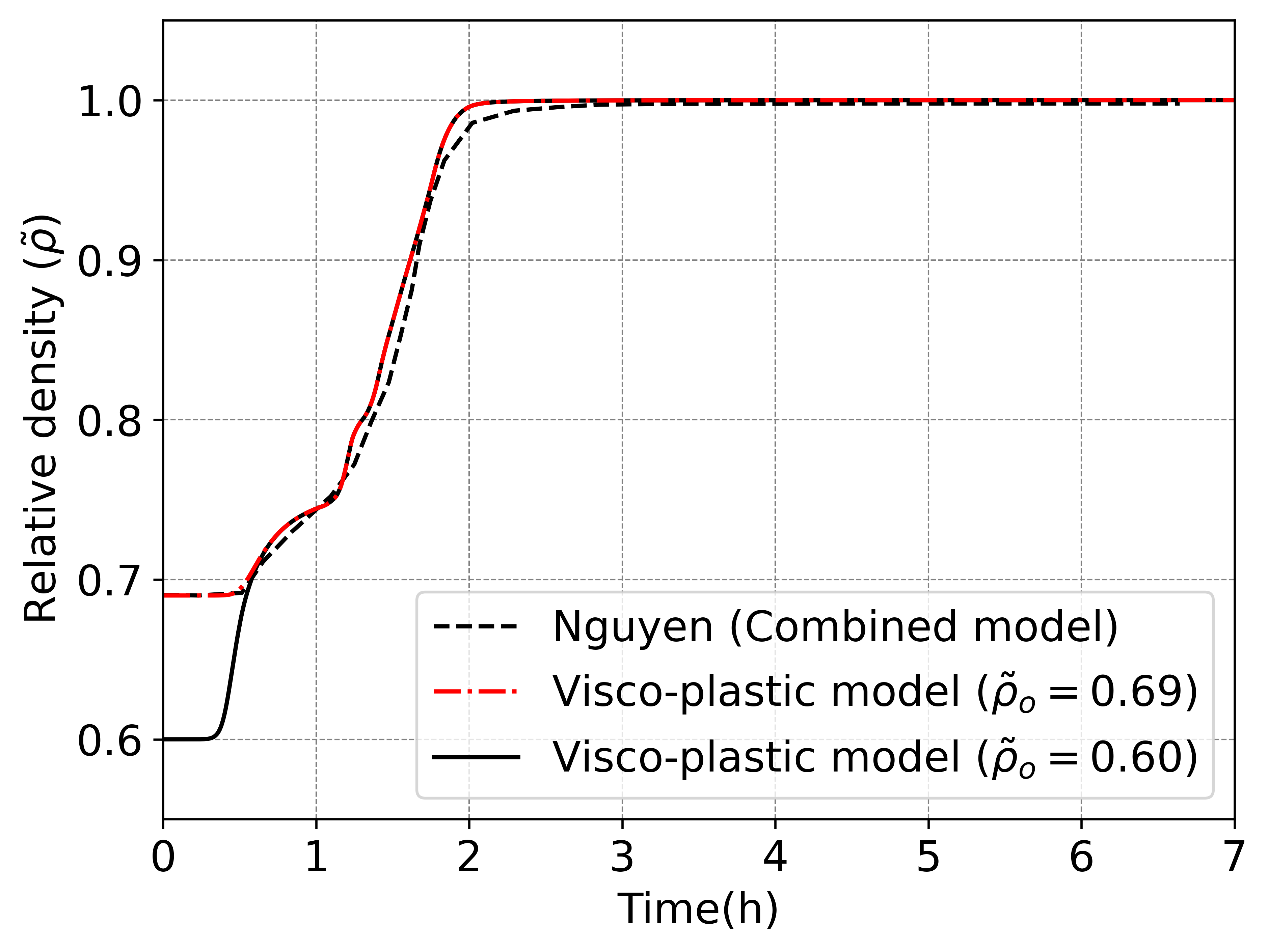}
        \caption{Visco-plastic model}
        \label{}
    \end{subfigure}
    \caption{The powder densification response of the calibrated \textbf{(a)} plastic model and \textbf{(a)} visco-plastic model with different initial relative densities ($\RD_o$) compared to the reference data (Nguyen \textit{et al.} \cite{vanthe2016numerical}).}
    \label{fig:plascompPlots}
\end{figure}

The initial relative density used for calibration is $\RD_o=0.69$ for both models (\fig \ref{fig:plascompPlots}) based on the Nguyen \textit{et al.} \cite{vanthe2016numerical} densification curve. However, in practical cases, much lower initial relative densities are encountered. Therefore, it is crucial for a PM-HIP model to account for different initial relative densities in the prediction of the densification behavior. A comparison of the response of plastic and visco-plastic models when subjected to lower initial relative densities ($\RD_o=0.60$) is shown in \fig \ref{fig:plascompPlots}.

It is observed that the plastic model's prediction (red dash-dot line) of the final relative density ($\RD_f$) is significantly reduced. This reduction can be explained by considering the expression of $\RD_f$ defined based on the volumetric plastic strain ($\epsilon^p_{vol}$) as,
\begin{equation}
    \RD_f = \RD_o e^{-\epsilon^p_{vol}} 
\end{equation}
In absence of any rate-dependent effects, the plastic (rate-independent) model undergoes a similar amount of $\epsilon^p_{vol}$ under the same loading conditions (pressure and temperature). Hence, with the same value of $\epsilon^p_{vol}$, but with a lower $\RD_o$ results in a lower $\RD_f$ prediction as shown in \fig \ref{fig:plascompPlots}a. This limitation of the plastic model is addressed by recalibrating the model with a lower $\RD_o$, while keeping the later part of the densification response unchanged, as shown in \fig \ref{fig:plascompPlots}a \cite{mayeur2025MCCP}. 

\begin{figure}[H] 
    \centering
    \begin{subfigure}[b]{0.36\textwidth}
        \centering
        \includegraphics[width=\textwidth]{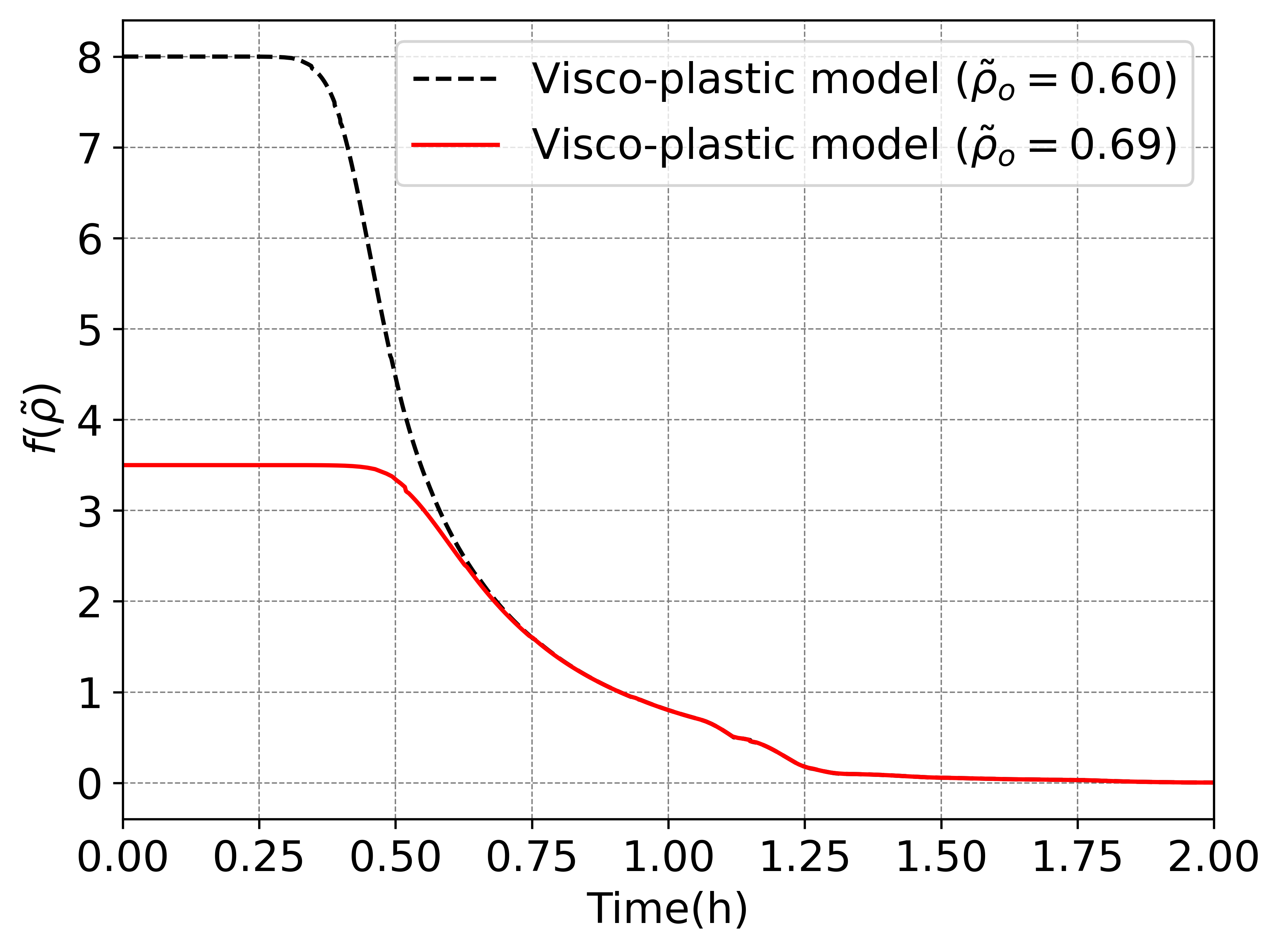}
        \caption{$f(\RD)$ vs Time}
        \label{}
    \end{subfigure}
    \hspace{1cm}
    \begin{subfigure}[b]{0.36\textwidth}
        \centering
        \includegraphics[width=\textwidth]{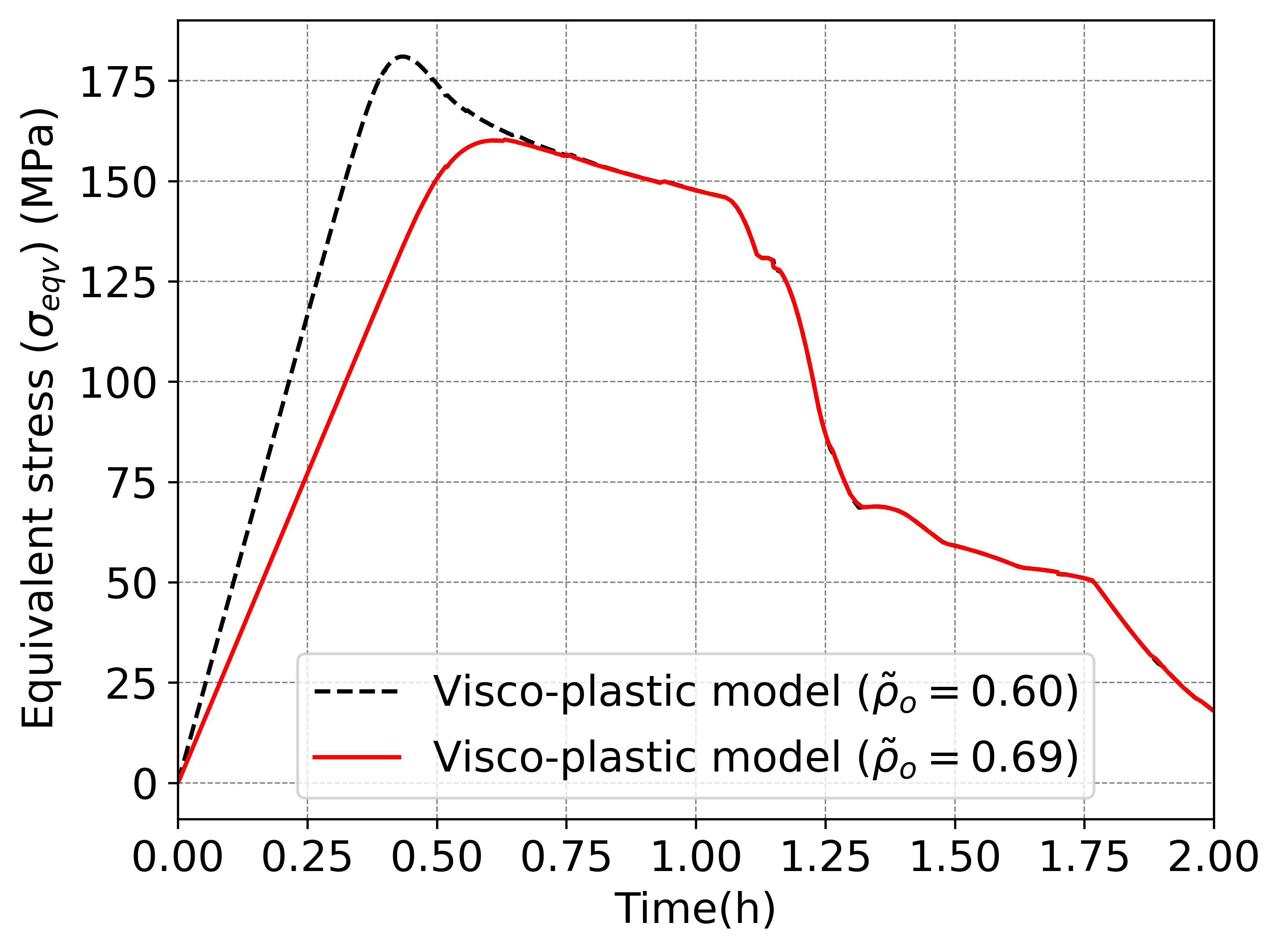}
        \caption{$\sigma_{eqv}$ vs Time}
        \label{}
    \end{subfigure}
    \caption{A comparison of the \textbf{(a)} $f(\RD)$ and \textbf{(b)} $\sigma_{eqv}$ during powder densification with different initial relative densities ($\RD_o$) in the visco-plastic model.}
    \label{fig:viscodesPlots}
\end{figure}

Contrastingly, the visco-plastic (rate-dependent) model readjusts the densification response based on the given $\RD_o$, as shown in \fig \ref{fig:plascompPlots}b. The densification in the visco-plastic model is driven by the volumetric visco-plastic strain ($\epsilon^{vp}_{vol}$) and is primarily responsible for the model's self adjustment to a lower $\RD_o$. This can be understood through the expression of the incremental swelling strain ($=\epsilon^{vp}_{vol}$) recalled below,
\begin{equation}\label{eq:AbqDeltaStainSw_lower}
    \Delta\Bar{\epsilon}^{sw} = \left( 3f(\RD)I_1 A(T) {\sigma}_{eqv}^{N(T)-1} \right) \Delta{t}
\end{equation}
With a lower $\RD_o$ ($=0.60$), all the quantities in Eq. (\ref{eq:AbqDeltaStainSw_lower}) have the same values (as with $\RD_o=0.69$) during the initial stages of densification except for $f(\RD)$ and $\sigma_{eqv}$ as shown in \fig \ref{fig:viscodesPlots}. Note that the $\sigma_{eqv}$ is affected due to its dependence on $f(\RD)$ (Eq. (\ref{eq:eqStressCr})). Hence, $f(\RD)$ is the only independent quantity with a different value at lower $\RD_o$. With a lower $\RD_o$ (higher $f(\RD)$), the visco-plastic model starts early densification and reaches the calibrated response. After reaching the calibrated response, the model with $\RD_o=0.60$ follows the same densification path as the calibrated model with $\RD_o=0.69$.

\subsection{Effects of stress shielding}
Stress shielding is the inability of the applied outside gas pressure to reach the inner sections of the powder due to can stiffness or a thicker cross-section of the component (more powder volume) \cite{lv2025review}. This causes non-uniform densification rates in different sections of the powder, contributing to final shape distortions. Therefore, accounting for stress shielding might improve the shape predictions, particularly for irregularly shaped components. 

\begin{figure}[h]
     \centering
     \includegraphics[width=\textwidth]{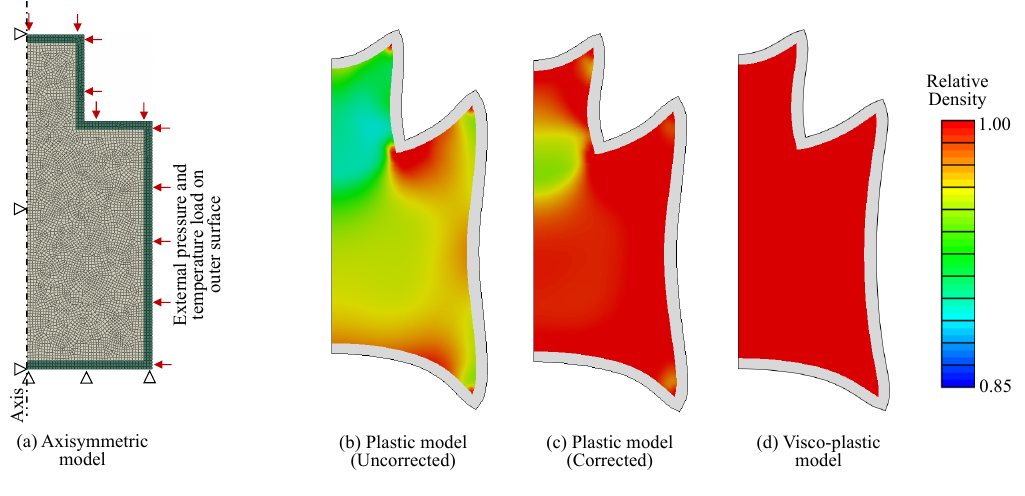}
     \caption{The predicted relative density ($\RD$) contour plots for the Nguyen \textit{et al.} \cite{van2017combined} geometry \textbf{(a)} using the plastic model [\textbf{(b)} - \textbf{(c)}] and the visco-plastic model \textbf{(d)}}
     \label{fig:StressShieldingContplots}
\end{figure}

\begin{figure}[H] 
    \centering
    \begin{subfigure}[b]{0.36\textwidth}
        \centering
        \includegraphics[width=\textwidth]{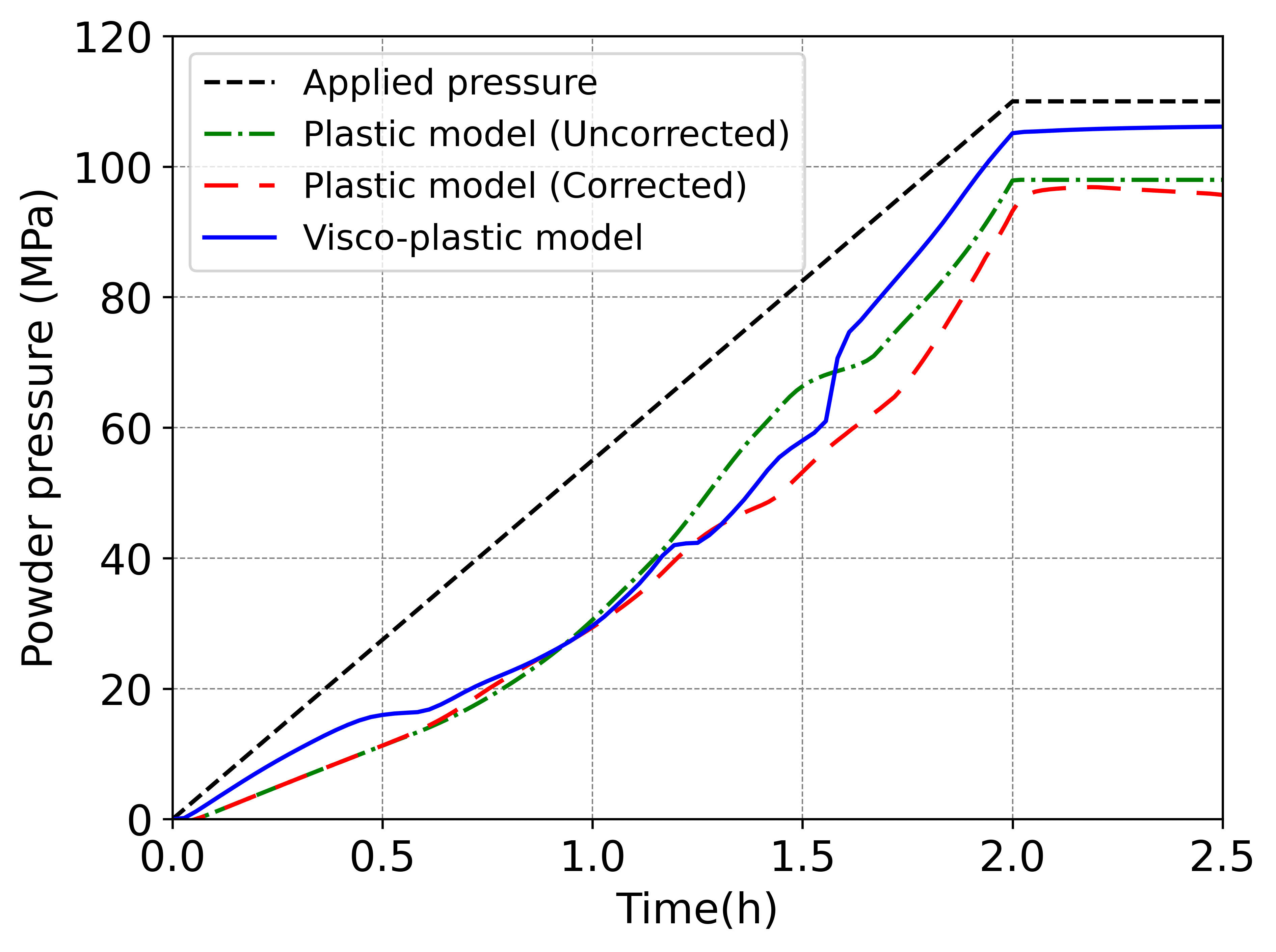}
        \caption{Mean powder pressure \textit{vs} Time}
        \label{}
    \end{subfigure}
    \hspace{1cm}
    \begin{subfigure}[b]{0.36\textwidth}
        \centering
        \includegraphics[width=\textwidth]{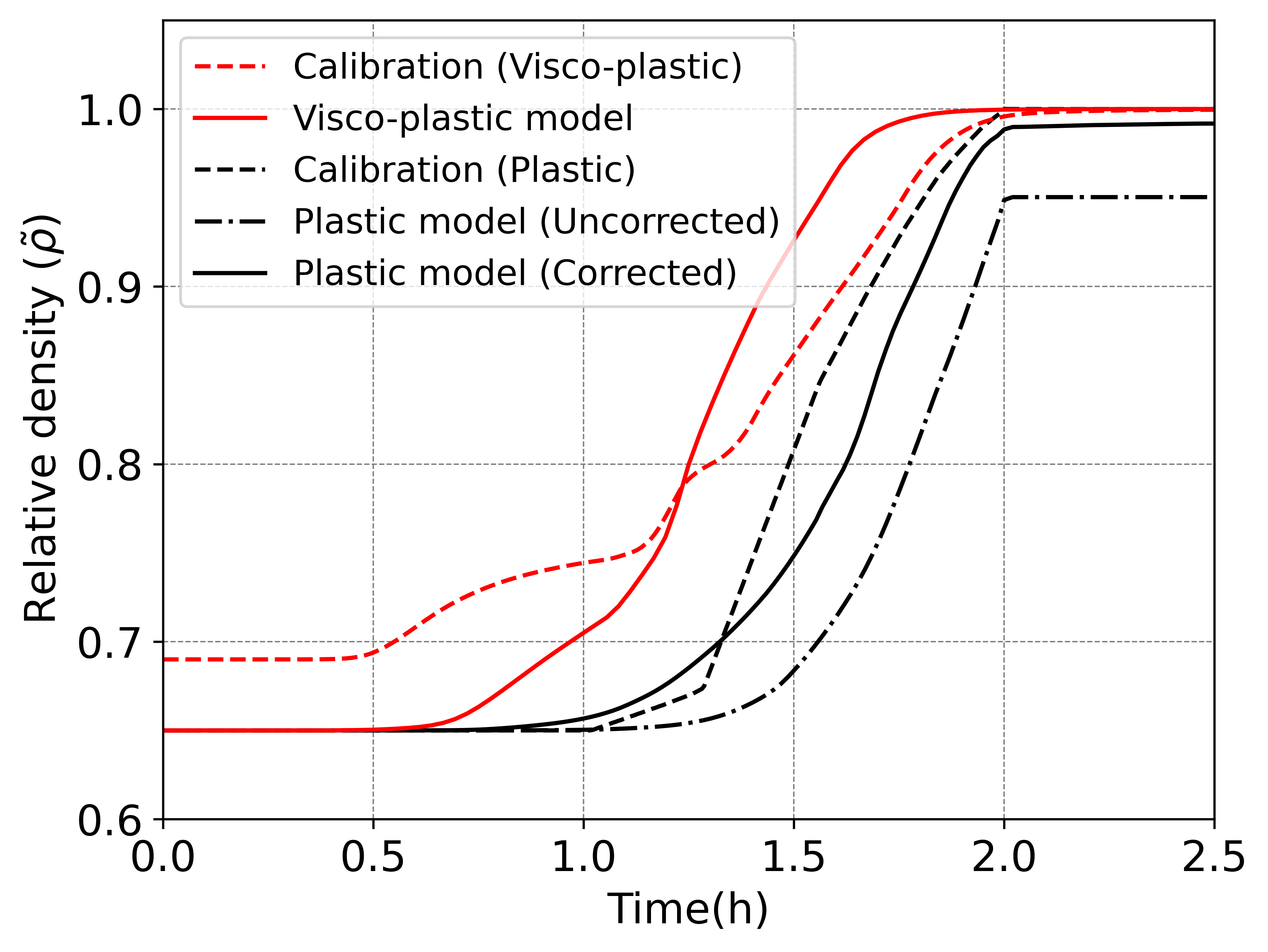}
        \caption{Mean $\RD$ \textit{vs} Time}
        \label{}
    \end{subfigure}
    \caption{A comparison of the \textbf{(a)} mean powder pressure and \textbf{(b)} mean relative density $(\RD)$ during densification in the plastic and visco-plastic model.}
    \label{fig:shieldingcompPlots}
\end{figure}

The effect of stress shielding on the plastic and visco-plastic models is studied through a reference problem in Nguyen \textit{et al.} \cite{van2017combined} shown in \fig \ref{fig:StressShieldingContplots}a. The applied pressure and temperature loads on the model are the same as those used for the single element simulations (\fig \ref{fig:singElem_HIPCyc}b). The predicted post-HIP relative density ($\RD$) contour plots obtained from the plastic and visco-plastic models are shown in \fig \ref{fig:StressShieldingContplots}b-d. \fig \ref{fig:StressShieldingContplots}b and \ref{fig:StressShieldingContplots}c depict the response of the plastic model, while the \fig \ref{fig:StressShieldingContplots}d shows the response of the visco-plastic model. Pronounced under-densification is observed in the plastic model (uncorrected) as compared to the visco-plastic model. This under-densification is caused by stress shielding, evident from the pressure \textit{vs} time plots shown in \fig \ref{fig:shieldingcompPlots}a. The reduced pressure experienced by the powder is observed to be equivalent in both the plastic model ($\approx$10-25 MPa) and visco-plastic model ($\approx$5-25 MPa). However, the overall stress shielding experienced during the entire ramp-up time (0-2 hours) remains higher in the plastic model than the visco-plastic model.

The effect of stress shielding can be understood in detail through the evolution of mean $\RD$ during the first two hours of densification as shown in \fig \ref{fig:shieldingcompPlots}b. For reference, the mean $\RD$ \textit{vs} time responses of the single element calibration model for both plastic and visco-plastic models are shown as dashed lines in \fig \ref{fig:shieldingcompPlots}b. Note that the plastic model calibrated with $\RD_o=0.60$ is used (black dashed line), while the visco-plastic model uses the calibration with $\RD_o=0.69$ (red dashed line). The mean $\RD$ \textit{vs} time response corresponding to the post-HIP contour plot in \fig \ref{fig:StressShieldingContplots}b is shown as a black dash-dot line in \fig \ref{fig:shieldingcompPlots}b labeled as plastic (uncorrected). This curve shows a significant lag in the $\RD$ evolution that leads to a large under-densification at the end of the ramp-up time (after 2 hours). This under-densification is addressed by applying a correction to the plastic model that involves repeated trials of this problem until the single element calibration response is approximately obtained (black solid line in \fig \ref{fig:shieldingcompPlots}b). The recalibrated (corrected) plastic model attains higher final $\RD$ as seen in the post-HIP $\RD$ contour and mean $\RD$ \textit{vs} time plots.

Contrastingly, the final $\RD$ predictions from visco-plastic model are not affected by lower pressures caused by stress shielding, as seen in the post-HIP contour plot (\fig \ref{fig:StressShieldingContplots}d). This is because the visco-plastic model attains full final relative density ($\RD_f=1.0$) well before the end of ramp-up time in the current problem, as seen in \fig \ref{fig:shieldingcompPlots}b (red solid line). The early densification can be attributed to less stress shielding experienced by the visco-plastic model after 1 hour mark, evident from the higher powder pressure in \fig \ref{fig:shieldingcompPlots}a. Moreover, the mean $\RD$ \textit{vs} time response of the visco-plastic can model is significantly different from the single element calibration response, due to the presence of additional effects of can geometry and size that were absent in the calibration model.

\subsection{Densification and distortion predictions}
In this subsection, the calibrated plastic and visco-plastic models are used for post-HIP relative density and shape predictions on different actual components that are also experimentally HIPed. These comparisons evaluate the prediction accuracy of the models under different conditions, such as different can size, can wall thickness and complex geometric shapes. 
\subsubsection{Cylindrical cans with different size and thickness}
\begin{figure}[H]
     \centering
     \includegraphics[width=0.75\textwidth]{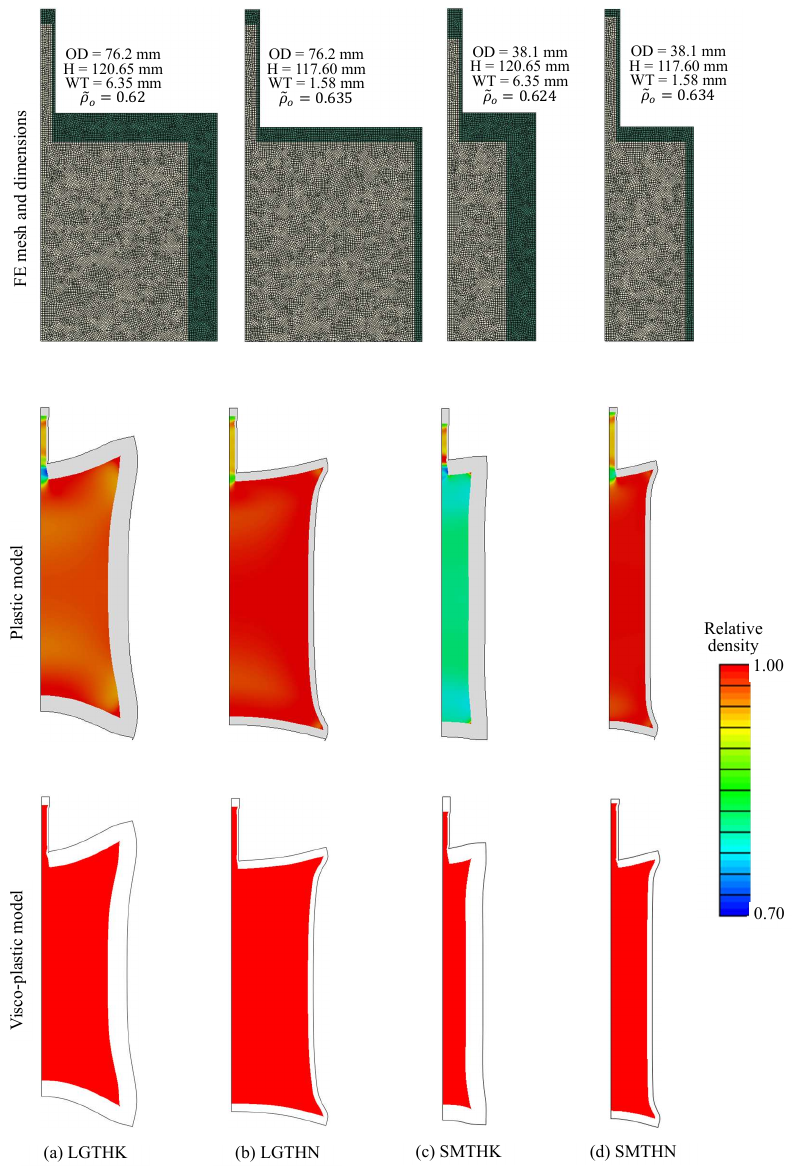}
     \caption{The FE mesh, dimensions (OD-Outer Diameter, H-Height, WT-Wall Thickness) and initial powder relative density ($\RD_o$) used for cylindrical can simulations: \textbf{(a)} Large thick can (LGTHK), \textbf{(b)} large thin can (LGTHN), \textbf{(c)} small thick can (SMTHK), \textbf{(d)} small thin can (SMTHN).}
     \label{fig:canMeshContPlots}
\end{figure}

A total of four cylindrical cans with different sizes and wall thicknesses are considered in this problem. The goal is to evaluate the effect of can size and can wall thickness on the model predictions. The dimensions of the cans, initial packing densities ($\RD_o$) and the finite element (FE) meshes used in simulations are shown in \fig \ref{fig:canMeshContPlots}. The experimental fabrication process, powder filling and HIPing process of the cans is described in Mayeur \textit{et al.} \cite{mayeur2025MCCP}. The can and powder material is SS316L, and a HIP cycle similar to that shown in \fig \ref{fig:singElem_HIPCyc}b is applied.

\begin{figure}[H] 
    \centering
    \begin{subfigure}[b]{0.36\textwidth}
        \centering
        \includegraphics[width=\textwidth]{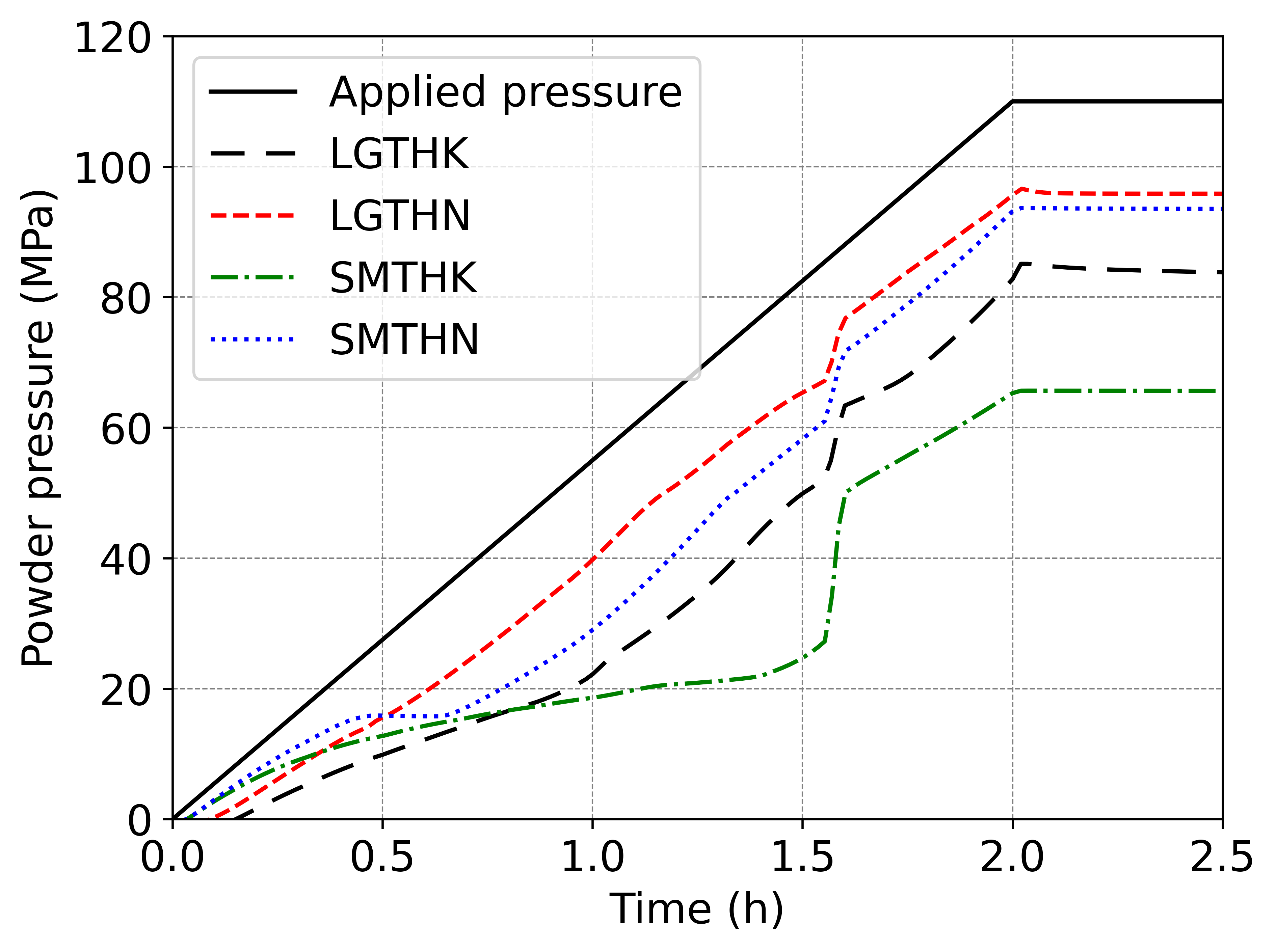}
        \caption{Mean powder pressure \textit{vs} Time (Plastic)}
        \label{}
    \end{subfigure}
    \hspace{1cm}
    \begin{subfigure}[b]{0.36\textwidth}
        \centering
        \includegraphics[width=\textwidth]{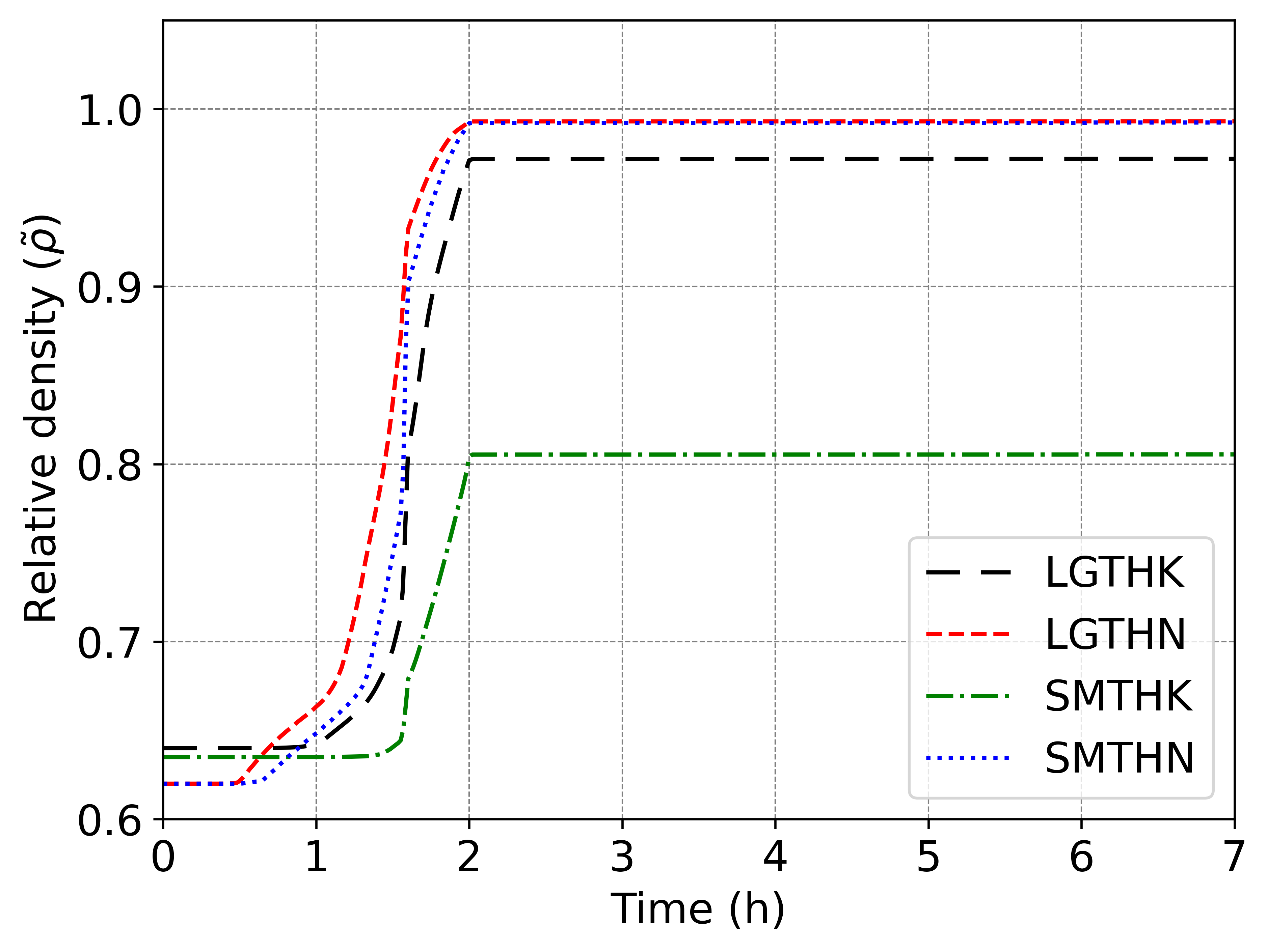}
        \caption{Mean $\RD$ \textit{vs} Time (Plastic)}
        \label{}
    \end{subfigure}
    \par\bigskip
    \begin{subfigure}[b]{0.36\textwidth}
        \centering
        \includegraphics[width=\textwidth]{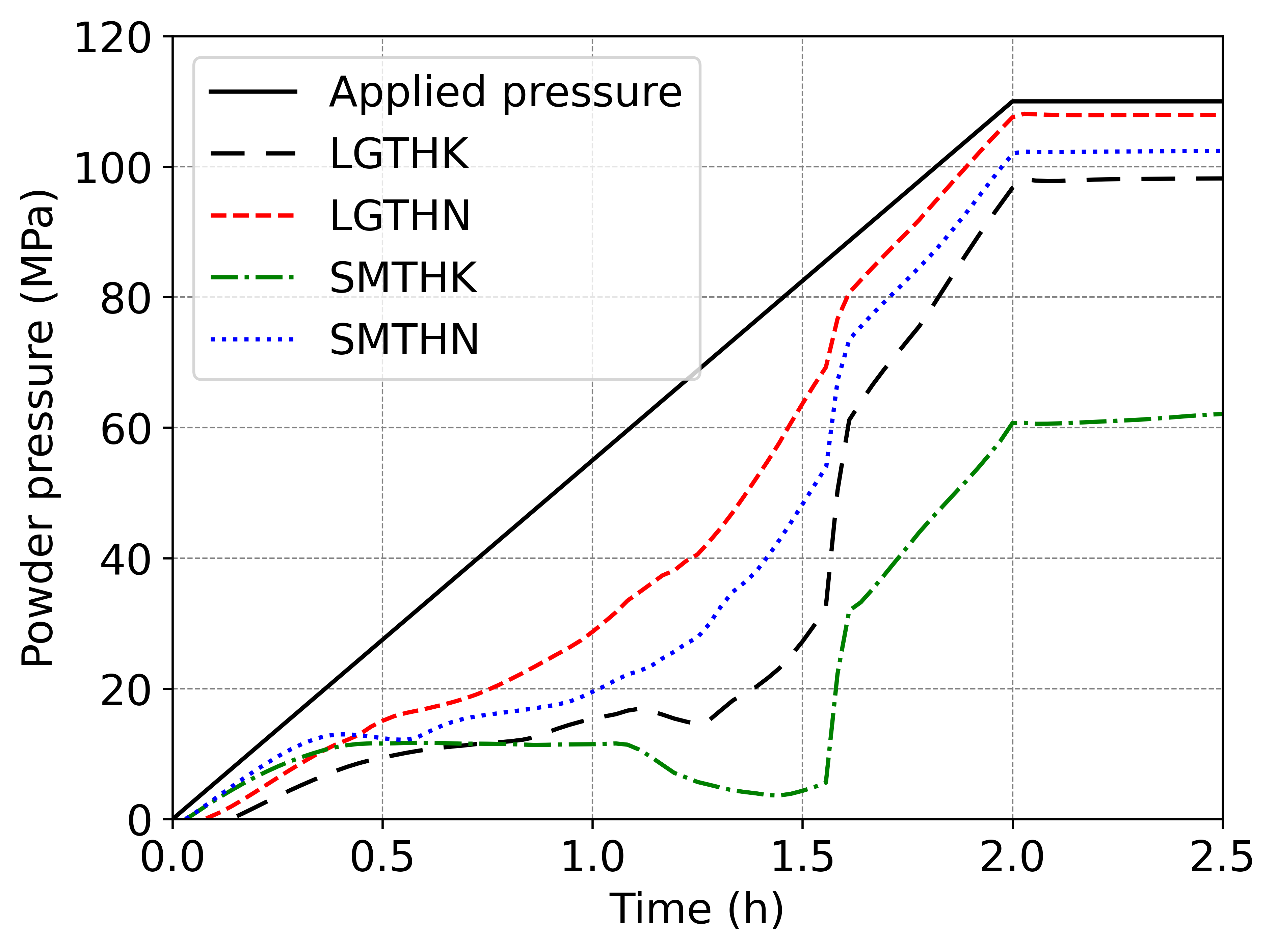}
        \caption{Mean powder pressure \textit{vs} Time (visco-plastic)}
        \label{}
    \end{subfigure}
    \hspace{1cm}
    \begin{subfigure}[b]{0.36\textwidth}
        \centering
        \includegraphics[width=\textwidth]{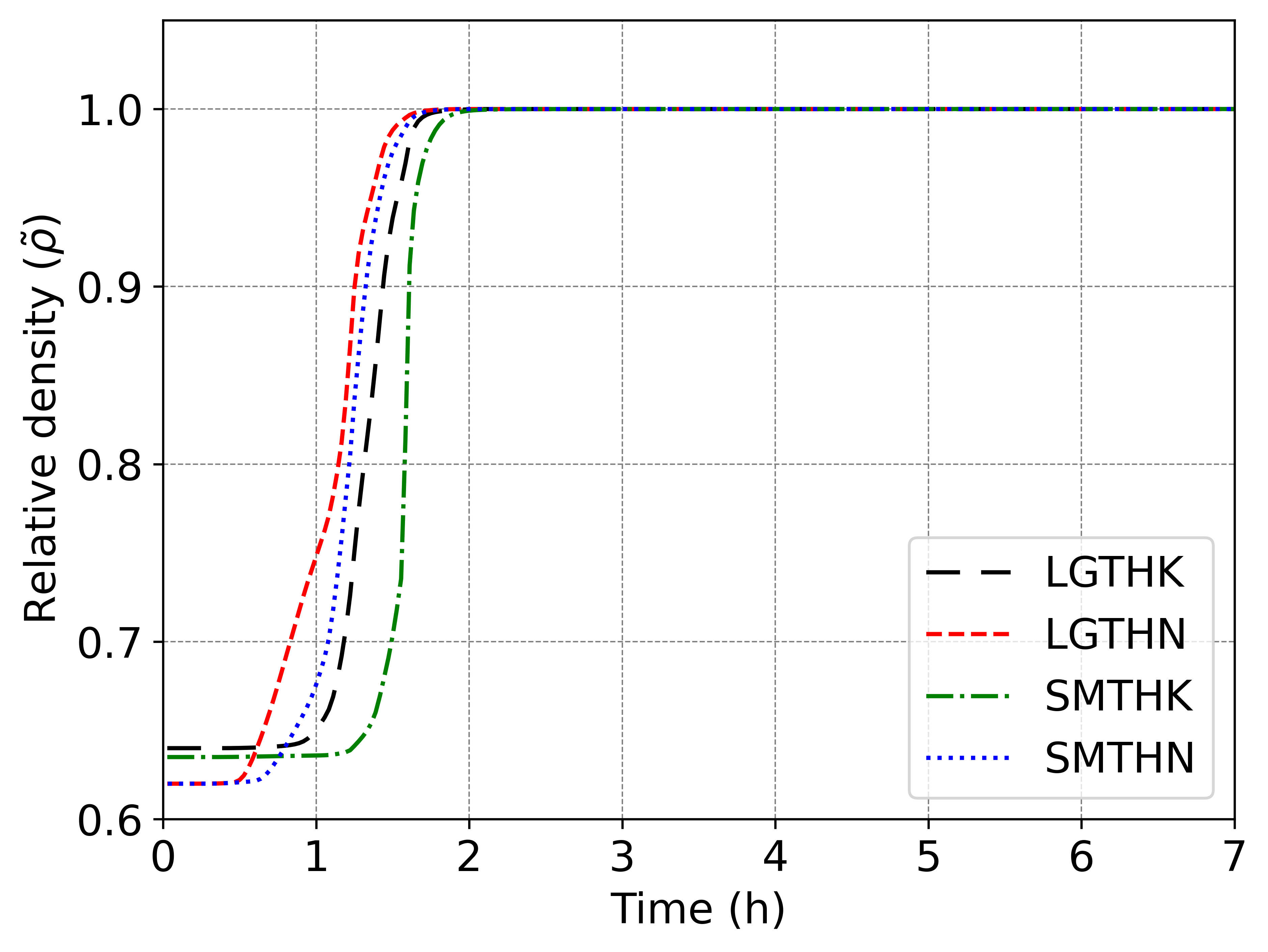}
        \caption{Mean $\RD$ \textit{vs} Time (visco-plastic)}
        \label{}
    \end{subfigure}
    \caption{A comparison of mean powder pressure and mean $\RD$ during the densification for different cylindrical cans using plastic (\textbf{a} - \textbf{b}) and visco-plastic (\textbf{c} - \textbf{d}) models.}
    \label{fig:viscocompPlots}
\end{figure}

The obtained relative density ($\RD$) predictions for the cans using plastic and visco-plastic models are shown in \fig \ref{fig:canMeshContPlots} for the four different cases named large thick can (LGTHK), large thin can (LGTHN), small thick can (SMTHK) and small thin can (SMTHN). Some under-densification ($\RD<1.0$) is observed in the plastic model cases, while all the cases in the visco-plastic model attain full relative density ($\RD=1.0$). Except for the small thick can (SMTHK), the under-densification in the plastic model is mostly confined to the fill stems with the rest of the can having homogenous relative density ($\RD\approx1.0$). The obsserved under-densification can be explained using the mean powder pressure and mean $\RD$ plots in \fig \ref{fig:viscocompPlots}. In the SMTHK can, under-densification of the plastic model can be attributed to higher stress shielding evident from the significantly lower powder pressure, as seen in \fig \ref{fig:viscocompPlots}a. Note that stress shielding is also present in the large thick (LGTHK) can, however, due to its larger size (and more powder volume), it experiences lower pressure reduction.

The densification in all the plastic model cases stops at the two-hour mark when the powder pressure stops increasing (\fig \ref{fig:viscocompPlots}b). This behavior is expected from a plastic (rate-independent) model, in which plastic deformation stops as soon as the stress state is on or within the yield surface. Contrastingly, all the cases of the visco-plastic model reach $\RD=1.0$ irrespective of the stress shielding. However, the case with higher powder pressure (LGTHN), i.e. less stress shielding, in the visco-plastic model attains $\RD=1.0$ earlier than the ones with lower powder pressure (SMTHK) (\fig \ref{fig:viscocompPlots}d). Therefore, the stress shielding tends to affect the transient densification behavior in a visco-plastic model rather than the final $\RD$ as in the plastic model.

The mean powder pressure plots for both the models in \fig \ref{fig:viscocompPlots}a and \fig \ref{fig:viscocompPlots}b show a peculiar drop in the rate of pressure increase starting at around 0.5 hours, followed by a pressure jump at around 1.5 hours. The decrease in the rate of pressure increase is found to be due to gradual thermal softening of the powder as the temperature increases. The pressure jump right after the thermal softening is found to be due to rapid compaction of powder that follows. The drops and jumps are more pronounced in the visco-plastic model than the plastic model because of slightly delayed start of densification in the visco-plastic model, evident from the mean $\RD$ \textit{vs} time plots. As a result of this delay, the powder has more time to soften in the visco-plastic model, hence, experiences lower powder pressures. Moreover, the cans with thicker walls experience more thermal softening of powder due to stress shielding, and hence have larger drops and jumps in the power pressure.

\begin{figure}[H]
    \centering
    \begin{subfigure}[b]{0.375\textwidth}
        \centering
        \includegraphics[width=\textwidth]{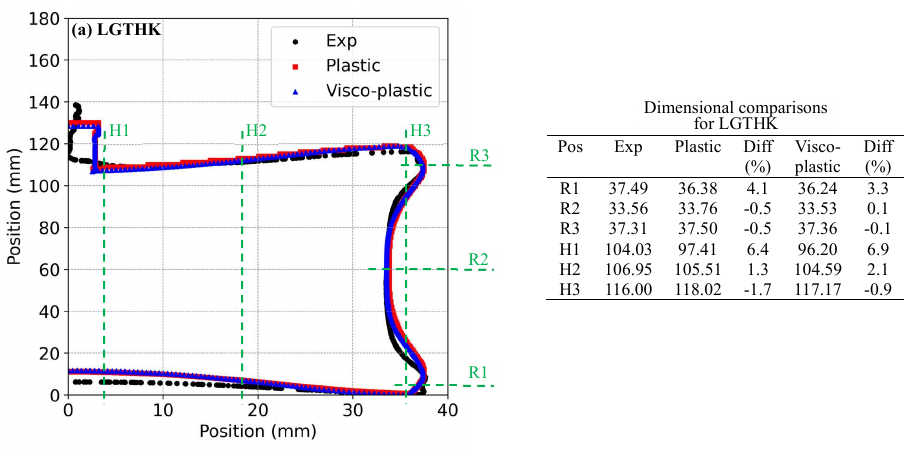}
        \label{}
    \end{subfigure}
    \hspace{1cm}
    \begin{subfigure}[b]{0.375\textwidth}
        \centering
        \includegraphics[width=\textwidth]{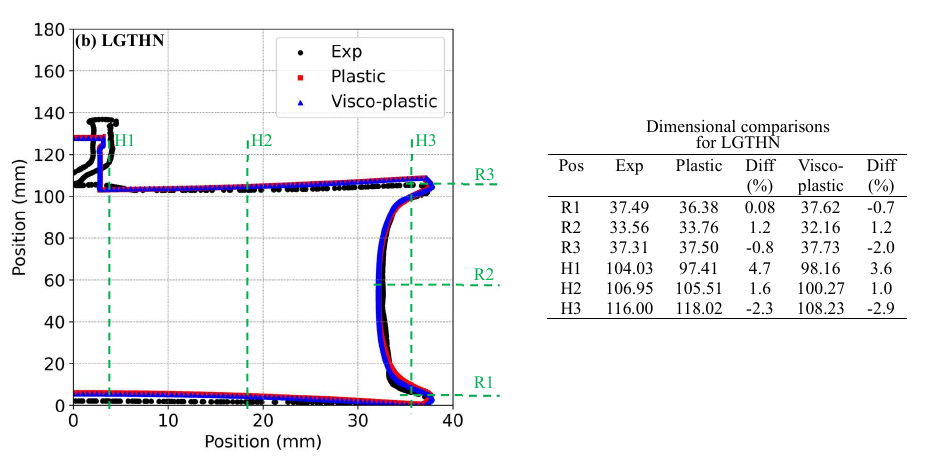}
        \label{}
    \end{subfigure}
    \vspace{0.25cm}
    \begin{subfigure}[b]{0.3\textwidth}
        \centering
        \includegraphics[width=\textwidth]{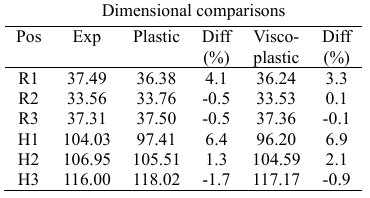}
        \caption{Large thick can (LGTHK)}
        \label{}
    \end{subfigure}
    \hspace{1.75cm}
    \begin{subfigure}[b]{0.3\textwidth}
        \centering
        \includegraphics[width=\textwidth]{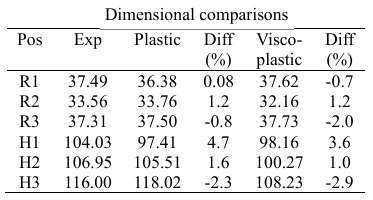}
        \caption{Large thin can (LGTHN)}
        \label{}
    \end{subfigure}
    \caption{The post-HIP shape comparisons for the large cylindrical can simulations obtained using plastic and visco-plastic models with experimental observations.}
    \label{fig:canShcompPlots1}
\end{figure}

\begin{figure}[H]
    \centering
    \begin{subfigure}[b]{0.375\textwidth}
        \centering
        \includegraphics[width=\textwidth]{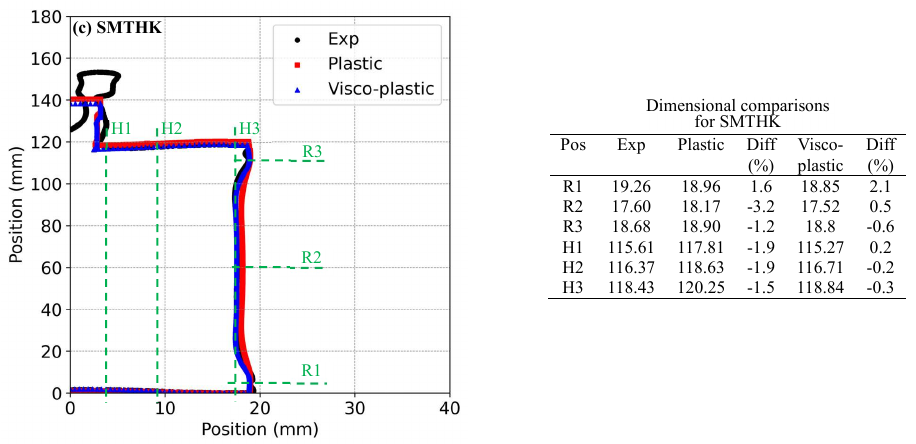}
        \label{}
    \end{subfigure}
    \hspace{1cm}
    \begin{subfigure}[b]{0.375\textwidth}
        \centering
        \includegraphics[width=\textwidth]{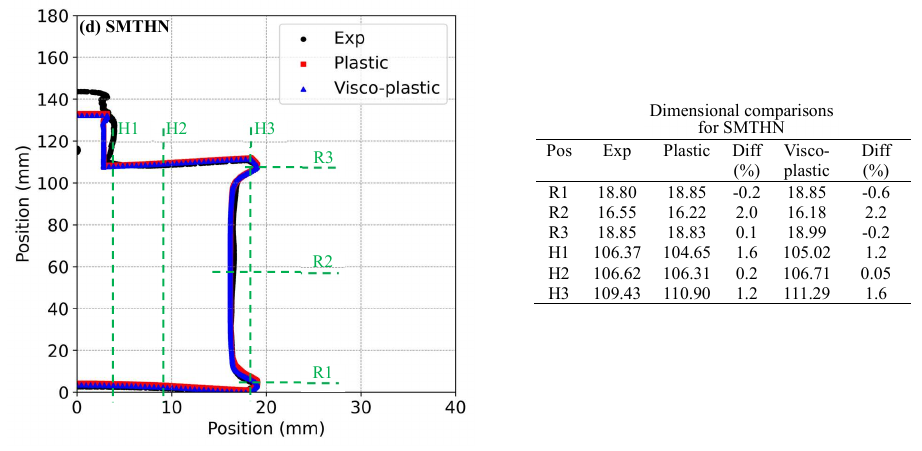}
        \label{}
    \end{subfigure}
    \vspace{0.25cm}
    \begin{subfigure}[b]{0.3\textwidth}
        \centering
        \includegraphics[width=\textwidth]{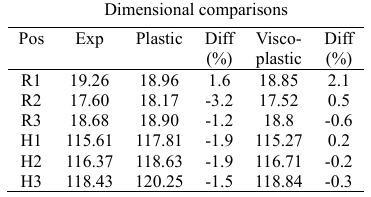}
        \caption{Small thick can (SMTHK)}
        \label{}
    \end{subfigure}
    \hspace{1.75cm}
    \begin{subfigure}[b]{0.3\textwidth}
        \centering
        \includegraphics[width=\textwidth]{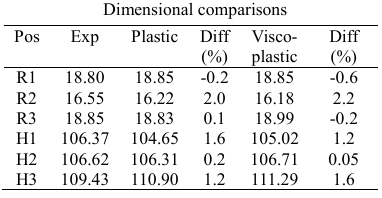}
        \caption{Small thin can (SMTHN)}
        \label{}
    \end{subfigure}
    \caption{The post-HIP shape comparisons for the small cylindrical can simulations obtained using plastic and visco-plastic models with experimental observations.}
    \label{fig:canShcompPlots2}
\end{figure}

Next, the shape comparisons are presented for all the four cases in \fig \ref{fig:canShcompPlots1} and \ref{fig:canShcompPlots2}. The shape predictions from the models are overlaid on the slices of 3D point clouds obtained from the experimental scans \cite{mayeur2025MCCP}. Visually, both the plastic and visco-plastic model predictions are almost equivalent and agree well with the experimental shape observations. The dimensional comparisons at various locations marked on the plots also indicate good agreement of the model predictions with experiments. The following observations are made from the plots and tabular comparisons in \fig \ref{fig:canShcompPlots1} and \ref{fig:canShcompPlots2},
\begin{itemize}
    \item The difference between the model predictions and experimental observations are mostly below 5\% except for the height measurement H1 of LGTHK can, where the deformation is over-estimated by both the models.
    \item The predictions of the radial dimensions are generally more accurate than the axial dimensions. The deviations in the predictions of axial dimensions from experiments is higher for larger cans.
    \item The end caps produce significant distortions in the final can shape that are more pronounced in larger cans with thinner walls.
\end{itemize}

\begin{figure}[H]
     \centering
     \includegraphics[width=0.9\textwidth]{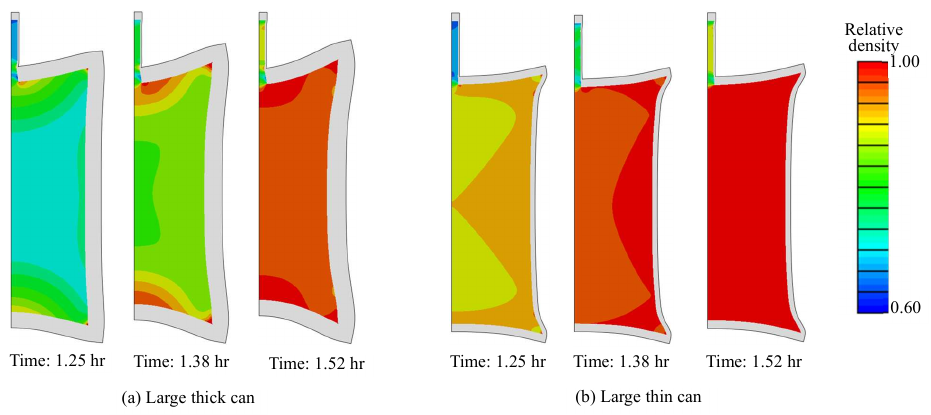}
     \caption{The relative density ($\RD$) contour plots for large thick (LGTHK) and large thin (LGTHN) cans during the densification at three intermediate time points.}
     \label{fig:caninterRD}
\end{figure}

The larger deviations in axial deformation predictions of the large cans can be explained through the $\RD$ contour plots at intermediate time points during densification shown in \fig \ref{fig:caninterRD} from the visco-plastic model. Differences in the $\RD$ evolution are evident from the plots of LGTHK and LGTHN cans. The densification rate in LGTHK is slower than the LGTHN which is also observed in the mean $\RD$ plots in \fig \ref{fig:viscocompPlots}d. Moreover, the densification front seen in the LGTHK can at 1.38 hours indicates that the densification starts from the end caps and moves inwards. As a result, at 1.38 hours, the end caps experience lower resistance to deformation and tend to deform inwards into the softer regions in LGTHK can ($\RD<1.0$). The lower resistance during the early stages of densification contributes to higher deviations observed in the final shapes of the larger cans. On the other hand, the smaller cans experience less end cap deformations, possibly due to higher can stiffness.

\subsubsection{Nozzle with stepped diameter}
\begin{figure}[H]
     \centering
     \includegraphics[width=0.9\textwidth]{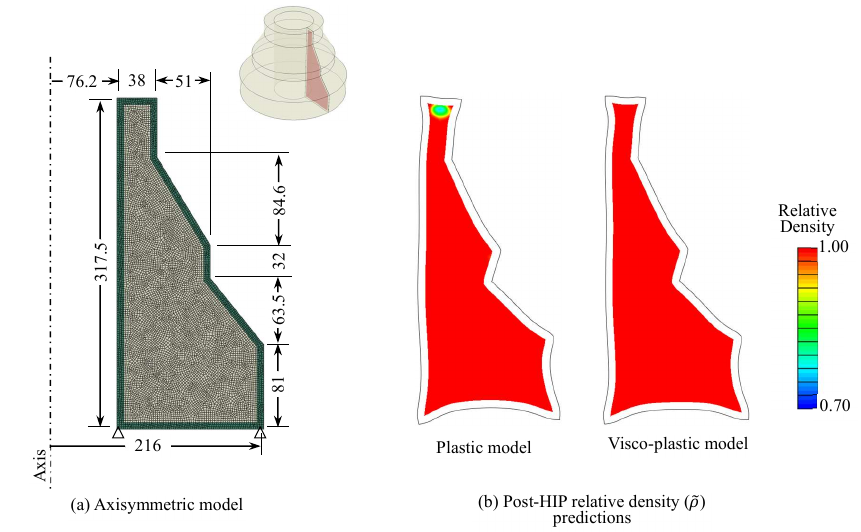}
     \caption{\textbf{(a)} The axisymmetric model with dimensions (in mm) used for nozzle simulation adopted from Sobhani \textit{et al.} \cite{sobhani2024thesis}. \textbf{(b)} A comparison of post-HIP relative density ($\RD$) predictions obtained from plastic and visco-plastic models.}
     \label{fig:EPRInozzle}
\end{figure}

In this example, a conical nozzle with stepped diameter is considered to evaluate the model predictions for a complex cylindrical geometry that has consecutive cylindrical and conical sections. The model description and experimental results are adopted from Sobhani \textit{et al.} \cite{sobhani2024thesis}. An axisymmetric model with dimensions and boundary conditions used for the simulation is depicted in \fig \ref{fig:EPRInozzle}a. The HIP cycle used for the simulation is similar to \fig \ref{fig:singElem_HIPCyc}b with the applied maximum pressure and temperature of 103 MPa and 1060 $^o$C respectively.

\begin{figure}[H]
     \centering
     \includegraphics[width=0.9\textwidth]{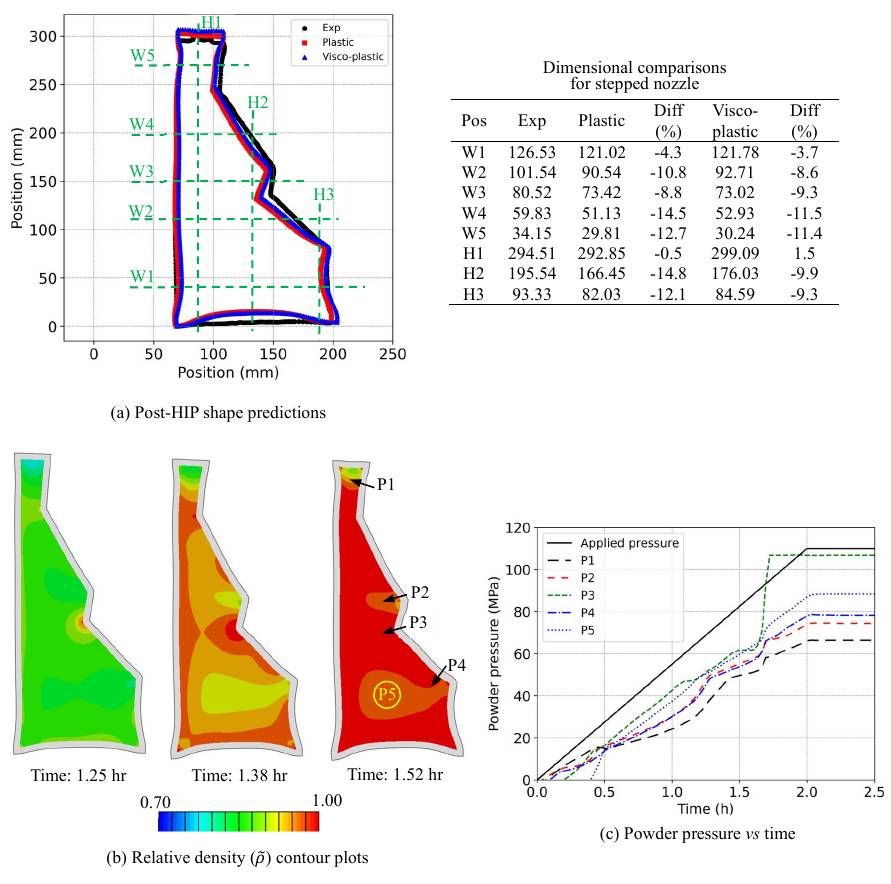}
     \caption{\textbf{(a)} A comparison of post-HIP shape predictions from the plastic and visco-plastic model with dimensional comparisons. \textbf{(b)} The relative density ($\RD$) contour plots at three intermediate times during densification. \textbf{(c)} The powder pressure \textit{vs} time at several regions inside the powder during the densification.}
     \label{fig:EPRINOZshapcomp}
\end{figure}

Final post-HIP relative density ($\RD$) predictions obtained for the nozzle using plastic and visco-plastic model are shown in \fig \ref{fig:EPRInozzle}b. The $\RD$ predictions show some under-densification in the plastic model, while the visco-plastic model undergoes full densification. The under-densification in the plastic model is confined to the top-most part where the powder cross-section is thinnest with two corners, causing higher stress shielding. \fig \ref{fig:EPRINOZshapcomp}a shows the post-HIP shape predicted by the plastic and visco-plastic models, which are found to be equivalent. However, notable deviations in the post-HIP shape predictions with respect to the experiments are observed for both the models at some locations. The table in \fig \ref{fig:EPRINOZshapcomp} indicates higher deviations in both the axial and radial dimensions, with the plastic model showing slightly higher deviations than the visco-plastic model. 

The higher values of dimensional deviations can be explained using the relative density ($\RD$) and the powder pressure plots shown in \fig \ref{fig:EPRINOZshapcomp}b and \ref{fig:EPRINOZshapcomp}c. From the $\RD$ plots (\fig \ref{fig:EPRINOZshapcomp}b), it is evident that the powder densification is non-uniform in the nozzle, with regions of faster and slower rates of densification. The non-uniform densification is caused by the nozzle's irregular geometry that creates regions of high and low stress shielding. The differences in stress shielding across different regions in the nozzle can be seen from the powder pressure variation during densification (\fig \ref{fig:EPRINOZshapcomp}c). In \fig \ref{fig:EPRINOZshapcomp}c, the powder pressures at five different regions (P1-P5) in the nozzle show that the region P3 experiences higher pressure than all the other marked regions. This is also observed in the $\RD$ plots (\fig \ref{fig:EPRINOZshapcomp}b) which show faster densification at region P3 than other regions. Moreover, the stress shielding is more pronounced in regions near the concave corners (P1, P2, P4), while the region near a convex corner (P3) experiences over pressure (opposite of stress shielding). The region P5 experiences lower pressure due to larger powder volume around it. This non-uniform pressure distribution leads to non-homogeneous densification of the powder. As a result of this non-homogeneous densification, the densified regions tend to move inwards into the under-densified softer regions during the early stages of densification, causing excessive distortions.

\subsubsection{T-valve}
\begin{figure}[H]
     \centering
     \includegraphics[width=0.6\textwidth]{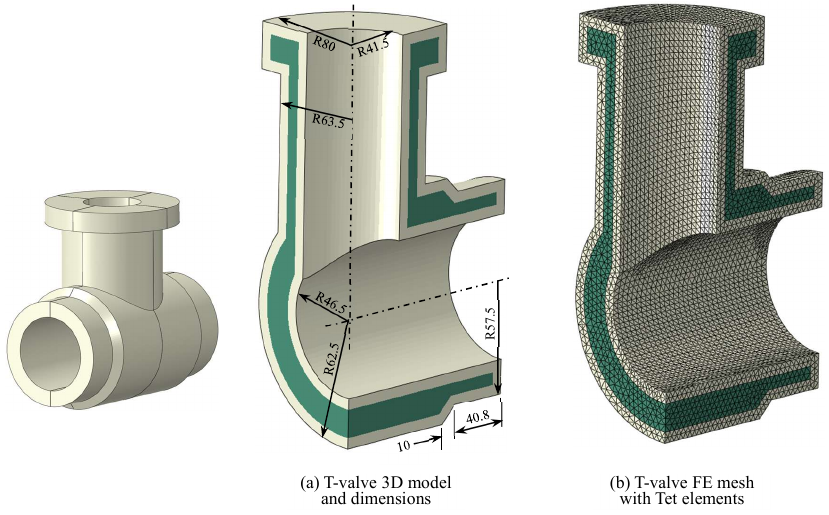}
     \caption{The 3D FE model used for the T-valve simulations with dimensions (in mm).}
     \label{fig:TVMesh}
\end{figure}

In this example, a three-dimensional (3D) geometry of a T-valve is chosen to evaluate the prediction accuracy of plastic and visco-plastic models for a more complex geometry that has more corners and curves. The quarter symmetric model of the T-valve with dimensions and finite element mesh is shown in \fig \ref{fig:TVMesh} which is adopted from Sobhani \textit{et al.} \cite{sobhani2024thesis}. The pressure and temperature applied for HIPing are 103 MPa and 1060 $^o$C respectively, following the same HIP cycle timing as in \fig \ref{fig:singElem_HIPCyc}b. The finite element mesh with approximately 99000 quadratic tetrahedral elements is used in the simulation.

\begin{figure}[H]
     \centering
     \includegraphics[width=0.95\textwidth]{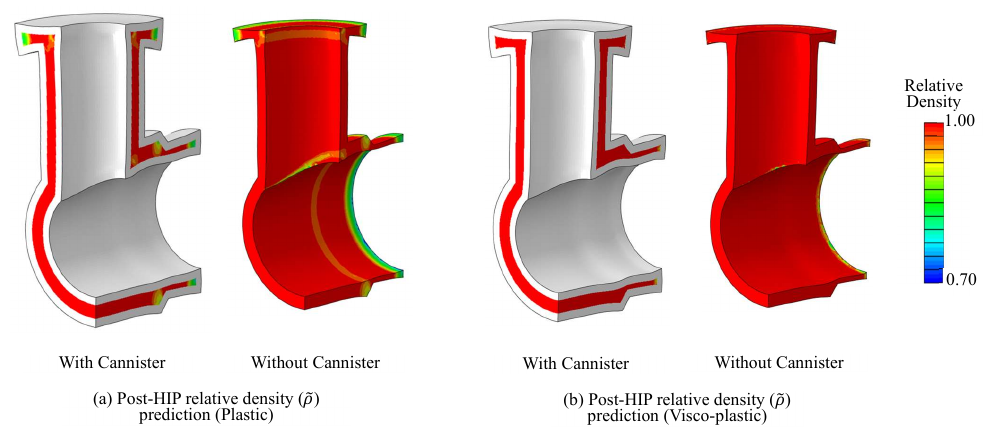}
     \caption{The post-HIP relative density ($\RD$) predictions for T-valve obtained using plastic and visco-plastic models.}
     \label{fig:TVRD}
\end{figure}

\begin{figure}[H]
     \centering
     \includegraphics[width=0.8\textwidth]{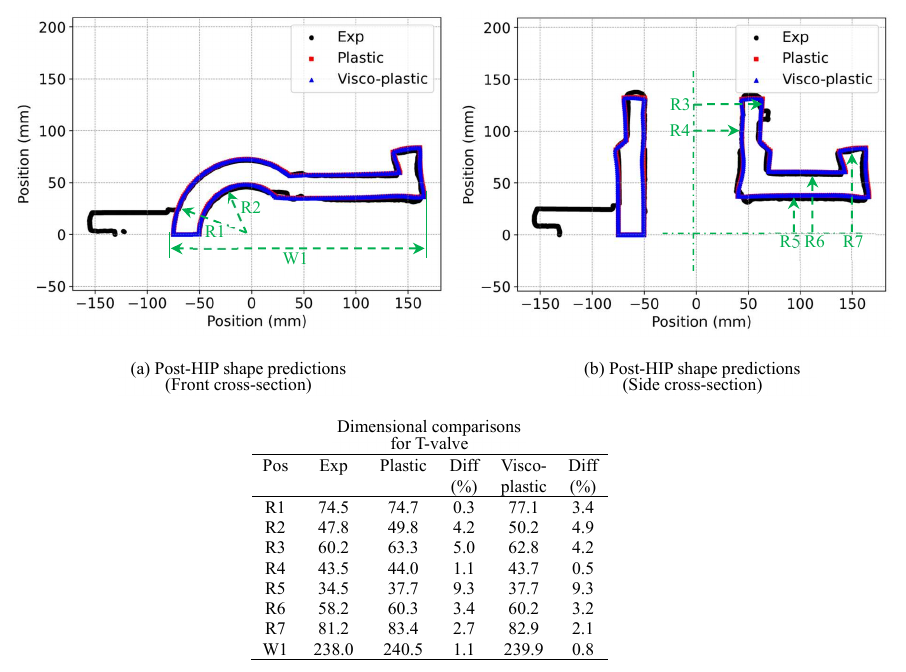}
     \caption{A comparison of post-HIP shape predictions of the T-valve using the plastic and visco-plastic models with the experimental observations.}
     \label{fig:TVShComp}
\end{figure}

The post-HIP relative density ($\RD$) plots for the T-valve obtained using plastic and visco-plastic model are shown in \fig \ref{fig:TVRD}. Some under-densification is observed in the plastic model at the extremities, while the visco-plastic model attains almost full densification. The regions of under-densification in both the models are near the concave corners. The post-HIP shape comparisons are shown in \fig \ref{fig:TVShComp} of two cross-section views of the T-valve. In both the views, most of the dimensional comparisons are within 5\% difference with respect to the experiments, indicating good agreement.

\section{Conclusion}
\label{Con}
In this work, a visco-plastic model is presented for simulating the powder metallurgy hot isostatic pressing (PM-HIP). A new modified calibration approach is employed for the visco-plastic model that utilizes less experimental data than existing approaches. With the new approach, the experimental data requirement is the same for both plastic and visco-plastic models, making their calibration cost equivalent. The visco-plastic model is implemented in a commercial FE package Abaqus using the CREEP subroutine and applied to simulate HIP on several different geometries. The following observations and insights are obtained from the results:

\begin{itemize}
    \item The visco-plastic and the plastic model calibrated with the same experimental data produce similar results in the majority of cases. However, differences are noticeable in the transient densification behavior, i.e., in relative density \textit{vs} time response.
    \item The visco-plastic model is insensitive to the slight variations in the HIP conditions such as pressure reduction due to stress shielding or the changes in the initial relative density, which may cause significant prediction errors in the plastic model.
    \item The plastic model, despite being sensitive to slight variations in the HIP conditions, produces accurate results after recalibration. In other words, the plastic model can also be made to handle these slight variations, albeit with some recalibration. 
    \item Both plastic and visco-plastic models are found capable of simulating HIP of complex geometries and predicting densification behavior accurately. The predictions of shape distortion from both the models agree reasonably well with the experimental observations.
\end{itemize}

It is already known and also demonstrated in this work that both the plastic and visco-plastic models have their own limitations and cannot fully represent the entire HIP process independently \cite{wikman2000combined,van2017combined}. For instance, the plastic model's accuracy suffers if the HIP conditions change, while the transient densification behavior in the visco-plastic model may deviate significantly from the reference behavior. Therefore, a unified model (plastic + visco-plastic) may be able to address these limitations and provide more accurate PM-HIP predictions.

\section*{Acknowledgement}
\noindent This work was supported by the Advanced Materials and Manufacturing Technologies Program of the U.S. Department of Energy’s Office of Nuclear Energy. The work was performed in partiality at the Oak Ridge National Laboratory’s Manufacturing Demonstration Facility, an Office of Energy Efficiency and Renewable Energy designated user facility. This manuscript has been authored by UT-Battelle, LLC, under contract DE-AC05-00OR22725 with the U.S. Department of Energy (DOE).

\section*{CRediT Author Statement}
\noindent \insertcreditsstatement

\appendix
\section{Derivatives}
\label{app_deriv}
For a computationally efficient implicit numerical integration, Abaqus requires the definition of several derivatives in the CREEP subroutine. These derivatives are named DECRA and DESWA representing derivatives of $\Delta\Bar{\epsilon}^{cr}$ and $\Delta\Bar{\epsilon}^{sw}$. The relevant derivatives for the present implementation are defined below,

\begin{equation}\label{eq:A1}
    \frac{\partial{\Delta\Bar{\epsilon}^{cr}}}{\partial{\Bar{\epsilon}^{sw}}} = A(T)\frac{\partial{c(\RD)}}{\partial{\Bar{\epsilon}^{sw}}}\sigma_{eqv}^{{N(T)-1}}q + A(T)c(\RD)(N-1)\sigma_{eqv}^{{N(T)-2}}\frac{\partial{\sigma_{eqv}}}{\partial{\Bar{\epsilon}^{sw}}}q
\end{equation}

\begin{equation}
    \frac{\partial{\Delta\Bar{\epsilon}^{cr}}}{\partial{p}} = A(T)\frac{\partial{c(\RD)}}{\partial{p}}\sigma_{eqv}^{{N(T)-1}}q + A(T)c(\RD)(N-1)\sigma_{eqv}^{{N(T)-2}}\frac{\partial{\sigma_{eqv}}}{\partial{p}}q
\end{equation}

\begin{equation}
    \frac{\partial{\Delta\Bar{\epsilon}^{cr}}}{\partial{q}} = A(T)c(\RD)\sigma_{eqv}^{{N(T)-1}} + A(T)c(\RD)(N-1)\sigma_{eqv}^{{N(T)-2}}\frac{\partial{\sigma_{eqv}}}{\partial{q}}q
\end{equation}

\begin{equation}
    \frac{\partial{\Delta\Bar{\epsilon}^{sw}}}{\partial{\Bar{\epsilon}^{sw}}} = -9A(T)\left[\frac{\partial{f(\RD)}}{\partial{\Bar{\epsilon}^{sw}}}\sigma_{eqv}^{{N(T)-1}}p
         + f(\RD)(N-1)\sigma_{eqv}^{{N(T)-2}}\frac{\partial{\sigma_{eqv}}}{\partial{\Bar{\epsilon}^{sw}}}p\right]
\end{equation}

\begin{equation}
    \frac{\partial{\Delta\Bar{\epsilon}^{sw}}}{\partial{p}} = -9A(T)\left[\frac{\partial{f(\RD)}}{\partial{p}}\sigma_{eqv}^{{N(T)-1}}p
         + f(\RD)(N-1)\sigma_{eqv}^{{N(T)-2}}\frac{\partial{\sigma_{eqv}}}{\partial{p}}p + f(\RD)\sigma_{eqv}^{{N(T)-1}}\right]
\end{equation}

\begin{equation}\label{eq:A6}
    \frac{\partial{\Delta\Bar{\epsilon}^{sw}}}{\partial{q}} = -9A(T)f(\RD)(N-1)\sigma_{eqv}^{{N(T)-2}}\frac{\partial{\sigma_{eqv}}}{\partial{q}}p
\end{equation}

\noindent All the quantities used in the Eqs. (\ref{eq:A1})-(\ref{eq:A6}) are defined in the constitutive model (Section \ref{contmodel}). The derivatives $\sigma_{eqv}$ appearing in Eqs. (\ref{eq:A1})-(\ref{eq:A6}) can be obtained using Eq. (\ref{eq:eqStressCr}). The derivatives of $c(\RD)$ and $f(\RD)$ with respect to $\RD$ are numerically computed from the calibration data using linear interpolation.

%
\bibliography{paper_bib.bib}
\end{document}